\documentclass[onefignum,onetabnum]{siamart251216} 

\usepackage{amssymb,mathtools,bm,stmaryrd,bbm,mathrsfs}
\SetSymbolFont{stmry}{bold}{U}{stmry}{m}{n}

\usepackage{tabularx,array,multirow,booktabs,diagbox}
\usepackage{float}
\usepackage[section]{placeins}
\usepackage{multicol}
\makeatletter
\let\caption@ORI@makecaption\@makecaption
\let\@makecaption\@undefined
\makeatother
\usepackage[compatibility=false]{caption}
\makeatletter
\let\@makecaption\caption@ORI@makecaption
\makeatother
\usepackage{subcaption}

\usepackage{graphicx}
\graphicspath{{fig/}}
\usepackage[export]{adjustbox}

\usepackage{xcolor}
\usepackage[normalem]{ulem}
\usepackage{comment}
\usepackage{enumerate}

\usepackage{siunitx}
\newsiamthm{problem}{Problem}
\newsiamremark{remark}{Remark}
\newsiamremark{example}{Example}

\newcommand{\vb}[1]{\bm{#1}}  %

\newcommand{\W}{\widetilde{W}}
\newcommand{\Winv}{\widetilde{W}^{-1}}

\headers{AL preconditioning for fictitious domain methods}{M.~Benzi, M.~Feder, L.~Heltai, and F.~Mugnaioni}

\title{A unified variational framework for augmented Lagrangian preconditioning of saddle point systems: applications to mixed-dimensional and fictitious domain methods\thanks{Submitted to the editors.\funding{MF and LH acknowledge support of the European Research Council (ERC) under the European Union's Horizon 2020 research and innovation programme (call HORIZON-EUROHPC-JU-2023-COE-03, grant agreement No. 101172493 "dealii-X").
The authors thank the Barcelona Supercomputing Center for the computational resources made available through the EuroHPC project EHPC-DEV-2026D02-275. The authors are members
of Gruppo Nazionale per il Calcolo Scientifico (GNCS), Istituto Nazionale di Alta Matematica (INdAM)}}}

\author{Michele Benzi\thanks{Scuola Normale Superiore, Piazza dei Cavalieri 7, 56126 Pisa, Italy (\email{michele.benzi@sns.it}, \email{federica.mugnaioni@sns.it}).}
\and Marco Feder\thanks{Department of Mathematics, University of Pisa, Largo B. Pontecorvo 5, 56126 Pisa, Italy (\email{marco.feder@dm.unipi.it}, \email{luca.heltai@unipi.it}).}
\and Luca Heltai\footnotemark[3]
\and Federica Mugnaioni\footnotemark[2]}

\begin{document}
\raggedbottom
\captionsetup[table]{skip=6pt}

\maketitle

\begin{abstract}
    We present a general framework to derive augmented Lagrangian preconditioners for saddle point systems arising from mixed finite element discretizations. Our approach is based on augmenting the variational formulation at the continuous level prior to discretization, providing a unified strategy for designing efficient preconditioners that are robust under mesh refinement for a broad class of problems leading to symmetric saddle point systems. The robustness of the proposed preconditioners is established through a detailed spectral analysis. This framework is applied to different model problems, focusing on mixed-dimensional and fictitious domain formulations involving Lagrange multipliers, with some configurations featuring large coefficient jumps.
    Extensive numerical experiments on different test cases and geometries are presented, validating the effectiveness and robustness of our preconditioners when used in conjunction with flexible GMRES.
\end{abstract}



\section{Introduction}\label{sec:intro}

The mathematical modeling of many physical and engineering problems leads to systems of partial differential equations subject to constraints. A classical approach to the treatment of such constraints consists in introducing Lagrange multipliers, which enforce the constraints in a weak sense within the variational formulation. After discretization with inf-sup stable elements, this procedure gives rise to algebraic saddle point systems of the form
\begin{equation}
    \label{eq:saddle-point}
    \begin{pmatrix} A & C^T \\ C & 0 \end{pmatrix}
    \begin{pmatrix} \mathbf{u} \\ \boldsymbol{\lambda} \end{pmatrix}
    =
    \begin{pmatrix} \mathbf{f} \\ \mathbf{g} \end{pmatrix},
\end{equation}
where $C$ is the discretization of the constraint operator, and $\boldsymbol{\lambda}$ denotes the vector of coefficients associated with the Lagrange multiplier. Besides enforcing the constraints, the multiplier often has a direct physical interpretation and its computation is therefore of independent interest.

Examples of fields where these large-scale saddle point problems naturally arise include computational fluid dynamics, constrained optimization, contact mechanics, optimal control, and many others. In recent years, mixed-dimensional models, in which the constraint operator $C$ couples problems posed on domains of different dimensionality, have gained increasing attention~\cite{trace1d2d,haznics,franco_zunino_dimola_numerical_2026}, as have fictitious domain methods, which give rise to a particularly challenging and relevant class of saddle point problems~\cite{Glowinski,BG_numermath,BOFFI20142151}. Rather than constructing a body-fitted mesh conforming to a geometrically complex or evolving physical domain, these methods embed the geometry of interest into a simpler \emph{background} domain that can be meshed independently, often by means of a structured grid. The interaction between the physical and background domains is then enforced weakly through a Lagrange multiplier field, resulting in
a saddle point formulation.
Depending on the support of the multiplier, one distinguishes between \emph{distributed} Lagrange multiplier (DLM) formulations, where the constraint is imposed through volumetric integrals, and \emph{boundary} Lagrange multiplier (BLM) formulations, where the multiplier is defined
only on the lower-dimensional interface between the background and immersed domains.
Fictitious domain methods based on distributed Lagrange multipliers have been successfully employed for elliptic interface problems~\cite{AURICCHIO201536,BOFFI20142151,BoffiGastaldi2024} as
a preliminary step to study fluid-structure interaction problems~\cite{BG_numermath}.
The resulting algebraic systems exhibit an intricate structure. In particular, the use of independent grids for the physical and background domains requires the assembly of coupling operators that encode interactions between basis functions
supported on non-matching and possibly overlapping meshes~\cite{BOFFI2022115650,boffi2023comparison}. As a consequence, the off-diagonal block $C$ in \eqref{eq:saddle-point} couples discrete spaces defined on different meshes, leading to saddle point systems whose analysis and implementation
are considerably more challenging than in standard boundary-fitted settings.

We cast all the problems discussed above into the following abstract framework, which we borrow from~\cite{boffi1997analysis}. For a Hilbert space $X$, we denote by $X'$ its dual and by
$\langle\cdot,\cdot\rangle_{X\times X'}$ the duality pairing, where the subscript indicates the spaces to which the first and second arguments belong. We denote by $(\cdot, \cdot)_X$ the inner product on $X$. Let $V$, $\Lambda$, and $H$ be real Hilbert spaces such that
\[
    \Lambda \subseteq H \equiv H' \subseteq \Lambda'
\]
is a Gelfand triple, with dense and continuous embeddings, where $H$ is identified
with its dual. This structure allows us to identify the duality pairing with the inner product of $H$ whenever
the second argument belongs to $H$, i.e., $\langle w,\mu\rangle_{\Lambda\times\Lambda'}=(w,\mu)_H$ for all $w\in\Lambda$ and
$\mu\in H$.

Let $a_0\colon V\times V\to\mathbb R$ be a continuous bilinear form and $L\colon V\to\Lambda$ a continuous linear operator. We consider the following classical
variational formulation: for $f\in V'$ and $g \in \Lambda$, find $(u,\lambda)\in V\times\Lambda'$ such that
\begin{equation}\label{eqn:problem_general}
    \begin{cases}
        a_0(u,v)+\langle Lv,\lambda\rangle_{\Lambda\times\Lambda'}
        =\langle f,v\rangle_{V'\times V}                 & \forall v\in V,         \\
        \langle Lu,\mu\rangle_{\Lambda\times\Lambda'}
        =\langle g,\mu \rangle_{\Lambda \times \Lambda'} & \forall \mu\in\Lambda'.
    \end{cases}
\end{equation}
It is well known that Problem~\eqref{eqn:problem_general} is well-posed if $a_0$ is coercive on the kernel $K=\{v\in V:Lv=0\}$ and $L$ is surjective~\cite{boffi2013mixed}.

The finite element approximation with $V_h\subset V$ and $\Lambda_h\subset H$ yields the discrete problem\footnote{In principle, a conforming approximation only requires $\Lambda_h \subset \Lambda'$. However, in most applications $H=L^2$, so
    any classical finite element space is naturally contained in $H$.}: find $(u_h,\lambda_h)\in V_h\times\Lambda_h$ such that
\begin{equation}\label{eqn:problem_general_discrete}
    \begin{cases}
        a_0(u_h,v_h)+\langle Lv_h,\lambda_h\rangle_{\Lambda\times \Lambda'}
        =\langle f,v_h\rangle_{V'\times V}                                                                      & \forall v_h\in V_h,        \\
        \langle Lu_h,\mu_h \rangle_{\Lambda \times \Lambda'}= \langle g, \mu_h\rangle_{\Lambda \times \Lambda'} & \forall \mu_h\in\Lambda_h.
    \end{cases}
\end{equation}
Classical results on saddle point problems~\cite{boffi2013mixed} show that problem~\eqref{eqn:problem_general_discrete} is uniquely solvable if there exist
mesh-independent positive constants $\alpha_h$ and $\beta_h$ such that
\begin{equation}\label{eqn:elker}
    a_0(v_h,v_h)\geq \alpha_h\|v_h\|_V^2
    \qquad \forall v_h\in K_h,
\end{equation}
\begin{equation}\label{eqn:inf-sup}
    \inf_{\mu_h\in\Lambda_h}\sup_{v_h\in V_h}
    \frac{\langle Lv_h,\mu_h\rangle_{\Lambda\times \Lambda'}}{\|v_h\|_V\,\|\mu_h\|_{\Lambda'}}
    =
    \inf_{\mu_h\in\Lambda_h}\sup_{v_h\in V_h}
    \frac{(Lv_h,\mu_h)_H}{\|v_h\|_V\,\|\mu_h\|_{\Lambda'}}
    \geq \beta_h,
\end{equation}where $K_h=\{v_h\in V_h:(Lv_h,\mu_h)_H=0\ \forall\mu_h\in\Lambda_h\}$.

\subsection{Motivation and main contributions}\label{sec:motivation}

Owing to the ubiquitous nature of large and sparse saddle point systems whose conditioning may deteriorate under mesh refinement and with variations in the physical and model parameters, a large variety of iterative solvers and preconditioning techniques have been developed in the literature~\cite{saadbook, ElmanSilvesterWathen,Benzi2005,mardal_winther}.
Among them, augmented Lagrangian (AL) techniques have emerged as a particularly powerful strategy. The basic idea consists in replacing the original system~\eqref{eq:saddle-point} with the equivalent augmented formulation
\begin{equation}\label{eqn:gen_aug_sys}
    \begin{pmatrix}
        A+\gamma C^T W^{-1}C & C^T \\
        C                    & 0
    \end{pmatrix}
    \begin{pmatrix}
        \mathbf{u} \\
        \boldsymbol{\lambda}
    \end{pmatrix}
    =
    \begin{pmatrix}
        \mathbf{f}+\gamma C^T W^{-1}\mathbf{g} \\
        \mathbf{g}
    \end{pmatrix},
\end{equation}
where $W$ is a suitable symmetric positive definite (SPD) matrix and $\gamma>0$ is a positive real number. An effective block triangular preconditioner is
\begin{equation}\label{eqn:AL_prec}
    \mathcal{P}_\gamma=
    \begin{pmatrix}
        A+\gamma C^T W^{-1}C & C^T                \\
        0                    & -\frac{1}{\gamma}W
    \end{pmatrix}.
\end{equation}

With an appropriate choice of the matrix $W$, increasing $\gamma$ clusters the spectrum and often leads to convergence rates that are independent of the mesh size and remarkably robust with respect to the physical parameters. Successful applications include linearizations of the stationary Navier–Stokes equations~\cite{ALprec,Farrell2019}, as well as multiphysics problems such as the coupled Stokes-Darcy problem~\cite{StokesDarcy} and magnetohydrodynamic problems~\cite{Laakmann2022}.
Even if the ideal preconditioner is replaced in practice by an inexact version, in which the diagonal blocks (or the action of their inverses) are replaced by suitable approximations, the large kernel of $C^T W^{-1}C$ makes the augmented block difficult to invert efficiently. Achieving optimal complexity therefore typically requires specialized, problem-dependent multigrid methods~\cite{ALprec,Farrell2019,AL_Oseen_Frank}, often based on the ideas introduced in~\cite{schoeberl_1999}.
Modified AL variants have also
been proposed~\cite{modALprec}, while more general algebraic strategies were tackled in~\cite{BenziFaccio2024}.
In the fictitious domain setting, we recently investigated AL preconditioners in~\cite{BENZI2026118522,benzi2026augmentedlagrangianpreconditionersfictitious}; related ideas were subsequently extended to mixed-dimensional beam-solid coupling in~\cite{FIRMBACH2027119425}. In these works, the augmentation is introduced algebraically by modifying the discrete system as in~\eqref{eqn:gen_aug_sys}. Although effective, this approach may require the explicit formation of $C^T W^{-1}C$, leading to a substantial computational cost and reduced sparsity of the resulting global system.

Motivated by these limitations and following the abstract theory of augmented Lagrangian formulations for mixed finite element schemes developed in~\cite{boffi1997analysis}, we directly augment the formulation~\eqref{eqn:problem_general} with $\gamma (Lu-g,Lv)_{\Lambda}$, associated with the constraint $Lu=g$ in $\Lambda$. Moreover, exploiting the continuous inclusion $\Lambda\subseteq H$ and incorporating an appropriate mesh-dependent scaling, we use $\gamma h^{-\alpha} (Lu_h-g,Lv_h)_H$ as the augmentation term in the discrete formulation, employing the more practical inner product on $H$ rather than that on $\Lambda$.
The resulting augmented bilinear form, $a_0(u_h, v_h)+\gamma h^{-\alpha}(Lu_h, Lv_h)_H$, can be directly assembled and avoids the explicit construction of the triple matrix product $C^T W^{-1}C$, preserving locality of assembly, reducing memory requirements, and facilitating the design of efficient solvers for the augmented subproblems.

This approach allows us to derive augmented formulations whose associated preconditioners admit a unified spectral analysis under the assumption that the bilinear form $a_0(\cdot, \cdot)$ is symmetric and positive semidefinite. We identify general conditions on the augmentation operator $W$ ensuring favorable spectral properties of the preconditioned systems. In particular, $W$ is the matrix representation of the scaled Riesz map associated with a suitable inner product on the multiplier space. %

However, unlike algebraic augmentation, modifying the variational formulation may alter the discrete solution. In some cases, as for Grad-Div stabilization with Taylor--Hood elements~\cite{HeisterRapinGradDiv, MOULIN2019718}, the resulting error can be shown to vanish as the mesh is refined. %
In such cases, the continuous augmentation asymptotically retains the properties of the algebraic augmentation while providing a more favorable framework for preconditioning. In other cases, however, the presence of a power of the mesh size in the formulation may make this error non-trivial to control.
A possible alternative is to use the projected counterpart $(\pi_H(Lu_h-g),\pi_HLv_h)_H$ rather than $(Lu_h-g,Lv_h)_H$, where $\pi_H$ is the $H$-orthogonal projector onto $\Lambda_h$, together with piecewise-constant multiplier spaces, when feasible~\cite{ALprec,Farrell2019}. Indeed, we prove that the projected realization coincides with the algebraic augmentation, and therefore does not alter the discrete solution. Moreover, when piecewise-constant multiplier spaces are used, the projected term can be implemented efficiently, avoiding the computational cost of explicitly forming $C^TW^{-1}C$. A concrete example is presented in Section~\ref{sec:num_dlm}.

We apply the resulting framework to the imposition of Dirichlet boundary conditions through Lagrange multipliers~\cite{BabuskaLM}, BLM and DLM fictitious domain formulations~\cite{benzi2026augmentedlagrangianpreconditionersfictitious,BENZI2026118522}, and immersed linear elasticity, thereby unifying and extending our earlier AL preconditioners developed for fictitious domain methods~\cite{BENZI2026118522,benzi2026augmentedlagrangianpreconditionersfictitious} within a common variational framework amenable to a much wider class of saddle point problems.

\subsection{Outline}

The remainder of the paper is organized as follows. Section~\ref{sec:notation} introduces the notation, while Section~\ref{sec:AL_formulations} presents the augmented formulation. The relationship between continuous and algebraic augmentation is discussed in Section~\ref{sec:algebraic_prop}, followed by the development and spectral analysis of the AL preconditioners in Section~\ref{sec:preconditioning}. Section~\ref{sec:numerical_tests} reports numerical experiments aimed at validating the theoretical results and assessing the practical performance of the proposed preconditioning strategy. Finally, Section~\ref{sec:conc} summarizes the main findings and outlines directions for future research.

\section{Notation}\label{sec:notation}
Given an open and bounded domain $D$, we denote by $L^2(D)$ the space of square-integrable functions on $D$, endowed with the norm $\|\cdot\|_{0,D}$ induced by the inner
product $(\cdot,\cdot)_{D}$. Sobolev spaces are denoted by $W^{s,p}(D)$, where $s\in \mathbb{R}$ refers to the differentiability index and $p \in [1,\infty]$ is the integrability exponent. In the case $p=2$, we use the notation $H^s(D)=W^{s,2}(D)$, with
associated norm $\|\cdot\|_{s,D}$ and seminorm $|\cdot|_{s,D}$. Dual spaces and duality pairings are denoted as in Section~\ref{sec:intro}.
The spectrum of a matrix is indicated by $\operatorname{Spec}(\cdot)$, while
vectors appear in boldface, e.g., $\mathbf{v}$. We will use the letter $\lambda$ for the Lagrange multiplier in the weak formulation. Finally, the letter $c$ will refer
to a constant independent of the discretization parameters, possibly taking different values from one formula to another.

\section{Augmented Lagrangian formulation}\label{sec:AL_formulations}

We start by recalling the abstract framework developed by Boffi \& Lovadina in~\cite{boffi1997analysis} for augmented Lagrangian formulations of mixed finite element schemes. Adding the term $\gamma(Lu-g,Lv)_{\Lambda}$, with $\gamma>0$, to the first equation of~\eqref{eqn:problem_general} does not change the continuous solution, since $Lu=g$ in $\Lambda$. The corresponding discrete scheme seeks $(u_h,\lambda_h)\in V_h\times\Lambda_h$ such that
\begin{equation}\label{eqn:problem_general_discrete_augmented_Lambda}
    \begin{cases}
        a_0(u_h,v_h) + \gamma(Lu_h-g,Lv_h)_{\Lambda} + \langle Lv_h,\lambda_h \rangle_{\Lambda \times \Lambda'} = \langle f,v_h\rangle_{V'\times V} & \forall v_h \in V_h,         \\
        \langle Lu_h,\mu_h \rangle_{\Lambda \times \Lambda'} = \langle g,\mu_h \rangle_{\Lambda \times \Lambda'}                                    & \forall \mu_h \in \Lambda_h.
    \end{cases}
\end{equation}
This is the natural augmentation, since the constraint residual $Lu_h-g$ is measured in the norm of the space in which $L$ takes its values. If $a_0$ is positive semidefinite (as will be assumed throughout), and the assumptions ensuring the well-posedness
of the original problem~\eqref{eqn:problem_general} hold, then for every $\gamma>0$ it can be shown that the
coercivity condition~\eqref{eqn:elker} is automatically satisfied for the bilinear form $a_0(\cdot,\cdot) + \gamma (L\> \cdot, L\> \cdot)_\Lambda$. Hence, only the discrete inf-sup condition~\eqref{eqn:inf-sup} remains to
be verified~\cite[Section~2]{boffi1997analysis}.

In practice, computing the inner product on $\Lambda$ may be impractical, for instance when $\Lambda=H^{\frac12}(\cdot)$. We therefore consider instead the inner product on $H$. This substitution, however, requires an appropriate mesh-dependent scaling.

Following the second augmented formulation of~\cite[Section~3]{boffi1997analysis}, we use $\gamma h^{-\alpha}(Lu_h-g,Lv_h)_H$ with $\alpha\geq 0$ left unspecified. The spectral analysis in Section~\ref{sec:preconditioning} will determine the choice of $\alpha$. The resulting scheme seeks $(u_h,\lambda_h) \in V_h \times \Lambda_h$, for $\gamma>0$, such that
\begin{equation}\label{eqn:problem_general_discrete_augmented}
    \begin{cases}
        a_0(u_h,v_h) + \gamma h^{-\alpha}(Lu_h -g,Lv_h)_{H} + \langle Lv_h,\lambda_h \rangle_{\Lambda \times \Lambda'}   =\langle f,v_h\rangle_{V'\times V}  \> \> \> \> \> \>\>                                                 & \forall v_h \in V_h          \\

        \langle Lu_h,\mu_h \rangle_{\Lambda \times \Lambda'}                                            = \langle g,\mu_h \rangle_{\Lambda \times \Lambda'} \> \> \>\> \> \>\> \> \>\> \> \>\> \> \>\> \> \>\> \> \>\> \> \>\>\> & \forall \mu_h \in \Lambda_h.
    \end{cases}
\end{equation}

Throughout this section, we only assume $a_0$ to be positive semidefinite on $V$, that
$L\colon V\to\Lambda$ is continuous and surjective, and that the discrete inf-sup
condition~\eqref{eqn:inf-sup} holds.

Since the $H$-norm used in the augmentation is weaker than the $\Lambda$-norm,
well-posedness of~\eqref{eqn:problem_general_discrete_augmented} additionally requires
the following coercivity condition~\cite{boffi1997analysis}: there exist constants $\xi,\delta>0$, independent
of $h$, such that
\begin{equation}\label{eqn:coercivity_augmented_eq}
    a_0(v,v) + \xi \|Lv\|_H^2 \geq \delta \|v\|_V^2
    \qquad \forall v \in V.
\end{equation}
This condition is application-dependent and must be verified on a case-by-case basis.

\subsection{Algebraic problem}
Given basis functions $\{ \phi_i\}_{i=1}^n$ and $\{ \psi_k\}_{k=1}^m$ such that $V_h \coloneqq \text{span} \{ \phi_i\}_{i=1}^n$ and $\Lambda_h \coloneqq \text{span} \{ \psi_k\}_{k=1}^m$, we have that Problem~\eqref{eqn:problem_general_discrete_augmented} can be written in matrix form as the following algebraic problem with unknowns $\mathbf{u}$ and $\vb{\lambda}$:
\begin{equation}
    \begin{pmatrix}\label{eqn:algebraic_form}
        \widetilde{A} & C^T \\
        C             & 0
    \end{pmatrix}
    \begin{pmatrix}
        \mathbf{u} \\
        \vb{\lambda}
    \end{pmatrix}=
    \begin{pmatrix}
        \widetilde{\mathbf{f}} \\
        \mathbf{g}
    \end{pmatrix},
\end{equation}
where
\begin{eqnarray*}
    \widetilde{A}_{ij} &\coloneqq&   \tilde{a}(\phi_j, \phi_i) \qquad \qquad \qquad \qquad \qquad i,j=1,\dots,n \\
    C_{k i} &\coloneqq& (L \phi_i,\psi_k )_H  \quad \> \qquad \qquad \qquad  \qquad \>\>\> i=1,\dots,n \quad k = 1,\dots, m  \\
    \widetilde{\mathbf{f}}_{i} &\coloneqq& \langle f, \phi_i\rangle_{V' \times V} + \gamma h^{-\alpha}(g, L\phi_i)_H    \qquad i=1,\dots,n \\
    \mathbf{g}_{k} &\coloneqq& \langle  g, \psi_k \rangle_{\Lambda \times \Lambda'} \qquad\qquad \qquad \qquad \quad \>\>\> k = 1,\dots, m,
\end{eqnarray*}and $\tilde{a}(\phi_i, \phi_j) \coloneqq a_0(\phi_i, \phi_j) + \gamma h^{-\alpha}(L\phi_i, L\phi_j)_H$ denotes the augmented bilinear form.

In the next section, we will take a closer look at the augmented term $(Lu_h,Lv_h)_H$, describing its finite dimensional representation in terms of known matrices.

\section{Connection between continuous and algebraic augmentation}\label{sec:algebraic_prop}

In this section, we show that the augmented term $(Lu_h,Lv_h)_{H}$ appearing in the augmented bilinear form $\tilde{a}(u_h, v_h) = a_0(u_h, v_h) + \gamma h^{-\alpha}(Lu_h, Lv_h)_H$, shares the same algebraic properties as the triple matrix product $C^T M^{-1}C$, where $M$
is the mass matrix defined on the multiplier space.
This observation generalizes the results of~\cite{HeisterRapinGradDiv}, where the specific case of Grad-Div stabilization for the Oseen problem is considered.

This naturally raises the question of whether considering the augmented formulation in~\eqref{eqn:problem_general_discrete_augmented} is equivalent to augmenting the discrete saddle point system as in~\eqref{eqn:gen_aug_sys} with $W=M$. As already observed in~\cite{HeisterRapinGradDiv} for Grad-Div stabilization, the two approaches are generally not equivalent. Nevertheless, we will show that
the term $(Lu_h,Lv_h)_{H}$ admits a decomposition into the sum of an algebraic term and a projection-type stabilization term, where the former coincides exactly with the entries of $C^T W^{-1}C$  when $W=M$.

We emphasize that, also in this section, we only assume $a_0(\cdot,\cdot)$ to be positive semidefinite, but not necessarily symmetric.

Let $\pi_H\colon H \rightarrow \Lambda_h$ be the $H$-orthogonal projector, i.e.
\begin{equation}\label{eqn:deforth}
    (w - \pi_H w,\chi_h)_H = 0 \quad \forall \chi_h \in \Lambda_h.
\end{equation}For $u_h, v_h \in V_h$, the term $(Lu_h-g,Lv_h)_H$ can be split into the following sum:
\begin{equation}\label{eqn:split}
    (Lu_h -g,Lv_h)_H = (\pi_H (Lu_h -g), Lv_h)_H + (\kappa(Lu_h -g),Lv_h)_H,
\end{equation}where $\kappa \coloneqq I- \pi_H$ is the so-called fluctuation operator. By the projection property of $\pi_H$, we have:
\begin{equation}\label{eqn:split2}
    (Lu_h-g,Lv_h)_H = (\pi_H (Lu_h - g), \pi_H (Lv_h))_H + (\kappa(Lu_h -g),\kappa(Lv_h))_H.
\end{equation} Following~\cite{HeisterRapinGradDiv}, we call the first summand \emph{algebraic term} and the second one \emph{stabilizing term}.
The following lemma shows that, after moving the term $(g,Lv_h)_H$ to the right-hand side of the first equation in~\eqref{eqn:problem_general_discrete_augmented}, the contribution $(\pi_H (Lu_h), \pi_H (Lv_h))_H$ can be represented as a product of known and computable matrices.

\begin{lemma}\label{lemma:algebraic_term}
    Let $\{\phi_i\}_{i=1}^n$ be a basis of $V_h$ and $\{\psi_k\}_{k=1}^m$ a basis of $\Lambda_h$. Let $C$ and $M$ be the following matrices, defined as:
    \[
        C_{k i} \coloneqq (L\phi_i, \psi_k)_H,
        \qquad
        M_{k l} \coloneqq (\psi_k, \psi_l)_H.
    \]
    Then,
    \[
        (\pi_H (L \phi_i),\; \pi_H (L \phi_j))_H=(C^T M^{-1} C)_{ij}.
    \]
\end{lemma}

\begin{proof}
    Let $w_{h_i}$ be the orthogonal projection of $L\phi_i$ onto $\Lambda_h$ with respect to the inner product $(\cdot, \cdot)_H$, that is $w_{h_i} \coloneqq\pi_H(L\phi_i) \in \Lambda_h.$
    Since $\{\psi_k\}_{k=1}^m$ is a basis of $\Lambda_h$, we may expand $w_{h_i}$ as $w_{h_i} = \sum_{k=1}^m \mathbf{m}_k^{(i)}\psi_k,$
    where $\mathbf{m}^{(i)}=(m_1^{(i)}, \ldots, m_m^{(i)}) \in \mathbb{R}^m$ is the corresponding coefficient column vector. By definition of
    orthogonal projection,
    \[
        (w_{h_i}, \chi_h)_H
        =
        (L\phi_i, \chi_h)_H
        \qquad \forall \chi_h \in \Lambda_h.
    \]

    Therefore, for each basis function $\psi_l$ we obtain
    \begin{equation}\label{eqn:system_MC}
        (w_{h_i}, \psi_l)_H
        =
        (L\phi_i, \psi_l)_H.
    \end{equation}
    Substituting the expansion of $w_{h_i}$ yields $\sum_{k=1}^m \mathbf{m}_k^{(i)} (\psi_k, \psi_l)_H = (L\phi_i, \psi_l)_H$. Moreover, using the definition of matrices $C$ and $M$ and letting $\mathbf{e}_i \in \mathbb{R}^n$ denote the $i$-th canonical basis vector, the
    system~\eqref{eqn:system_MC} can be written in matrix form as $M\, \mathbf{m}^{(i)} = C \mathbf{e}_i.$ We now compute the $(i,j)$-entry of $C^T M^{-1} C$:
    \[
        \mathbf{e}_i^T C^T M^{-1} C \mathbf{e}_j
        \!=\!(M^{-1} C \mathbf{e}_i)^T M (M^{-1} C \mathbf{e}_j)\!=\!\bigl(\mathbf{m}^{(i)}\bigr)^T M \mathbf{m}^{(j)}
        \!=\!(w_{h_i}, w_{h_j})_H,
    \]
    where the last equality follows from expanding both $w_{h_i}$ and $w_{h_j}$ in the basis $\{\psi_k\}_{k=1}^m$. Finally, recalling that $w_{h_i} = \pi_H (L\phi_i)$, we conclude that
    \[
        (C^T M^{-1} C)_{ij}
        =
        \bigl(\pi_H (L \phi_i),\; \pi_H (L\phi_j)\bigr)_H.
    \]
\end{proof}

The following lemma concerns instead the \textit{stabilizing term}.

\begin{lemma}\label{lemma:stab_term}
    Let $(u_h, \lambda_h)$ be a solution of the augmented Lagrangian formulation~\eqref{eqn:problem_general_discrete_augmented}, and let $(u, \lambda)$ solve the original
    problem~\eqref{eqn:problem_general}. We define the augmentation residual, namely, the term by which~\eqref{eqn:problem_general_discrete_augmented} differs from~\eqref{eqn:problem_general_discrete}, as
    \[
        r(v_h) \coloneqq \gamma h^{-\alpha} (L u_h- g, L v_h)_H, \qquad \forall v_h \in V_h.
    \]Then, the following bound holds:
    \begin{equation}\label{eqn:decay_stabilizing_term}
        |r(v_h)| \le \gamma h^{-\alpha}  \|\kappa (L (u_h- u))\|_H \|\kappa(L v_h)\|_H.
    \end{equation}
    In particular, if $L(V_h) \subset \Lambda_h$ and $g \in \Lambda_h$, then $r(v_h) \equiv 0$.
\end{lemma}
\vspace{0.15cm}
\begin{proof}
    The second equation of the augmented problem~\eqref{eqn:problem_general_discrete_augmented} gives
    \[(Lu_h-g,\mu_h)_H=0,
        \qquad \forall \mu_h\in\Lambda_h,\]
    since $(Lu_h -g) \in\Lambda\subseteq H$ and $\mu_h\in\Lambda_h\subset H$, so that the duality pairing between $\Lambda$ and $\Lambda'$ can be identified with the scalar product in $H$. Using the definition of $H$-orthogonal projection onto $\Lambda_h$~\eqref{eqn:deforth}, the relation above implies $\pi_H(Lu_h-g)=0.$
    Then, using the decomposition in~\eqref{eqn:split2}, we obtain
    $$
        r(v_h)
        =
        \gamma h^{-\alpha}
        \bigl(\kappa(Lu_h-g),\kappa(Lv_h)\bigr)_H,
    $$ which contains the stabilizing term introduced in~\eqref{eqn:split2}.
    Finally, using Cauchy-Schwarz inequality and the fact that $Lu_h-g=L(u_h -u)$ we obtain the desired bound
    \begin{align*}
        |r(v_h)| & = \gamma h^{-\alpha} |(\kappa (L u_h- g), \kappa(L v_h))_H|                                                                               \\
                 & \le \gamma h^{-\alpha} \|\kappa (L u_h- g)\|_H \|\kappa(L v_h)\|_H   = \gamma h^{-\alpha}  \|\kappa (L (u_h- u))\|_H \|\kappa(L v_h)\|_H.
    \end{align*}
\end{proof}

\begin{remark}\label{remark:stab_term_vanishes}
    The above bound is clearly problem-dependent. For instance, in the context of Grad-Div stabilization, where $L$ is the divergence operator and $g=0$, it was shown in~\cite{HeisterRapinGradDiv} that, for Taylor--Hood elements, $\|\kappa(\operatorname{div}(u_h - u))\|_{0,\Omega} \to 0$ as $h \to 0$.
    Whenever this occurs, we can conclude that the continuous augmentation
    asymptotically retains the same properties as the algebraic one, while avoiding the computational drawbacks associated with the explicit
    formation of $C^T M^{-1}C$. In Section~\ref{subsec:weak_bcs}, we will show another illustrative example for which the stabilizing term can
    be proved to vanish as the mesh is refined. A further possibility is to choose, if feasible, discrete spaces satisfying $L(V_h) \subset \Lambda_h$. This is the case for the classic Stokes problem discretized with divergence-free pairs, such as Scott--Vogelius elements.
\end{remark}

In the following section, we introduce the proposed AL preconditioner and provide a detailed spectral analysis that accounts for the effects of the stabilizing term arising from the continuous augmentation.

We define the symmetric positive semidefinite matrix $E$ associated with the bilinear form $(\kappa(Lu_h),\kappa(Lv_h))_H$, i.e.
\begin{equation}\label{eqn:matrix_E}
    \mathbf{v}^T E \mathbf{v}
    =
    \|\kappa (L v_h)\|_H^2 \qquad \forall v_h \in V_h.
\end{equation}

\section{Preconditioning and spectral analysis}\label{sec:preconditioning}

Using the decomposition~\eqref{eqn:split2} and Lemma~\ref{lemma:algebraic_term}, the matrix $\widetilde{A}$ in~\eqref{eqn:algebraic_form} can be decomposed as $\widetilde{A} = \hat{A} + \gamma h^{-\alpha} C^T M^{-1}C$, where $\hat{A} \coloneqq A+\gamma h^{-\alpha}E$ and $E$ is the matrix
in~\eqref{eqn:matrix_E}. We denote the matrix in system~\eqref{eqn:algebraic_form} by
\begin{equation}\label{eqn:aug_system_matrix}
    \widetilde{\mathcal{A}}_\gamma =
    \begin{pmatrix}
        \widetilde{A} & C^T \\
        C             & 0
    \end{pmatrix}.
\end{equation}
By analogy with the algebraic augmentation, it is natural to consider a corresponding preconditioner of the form
\begin{equation}\label{eqn:candidate_prec}
    \widetilde{\mathcal{P}}_\gamma =
    \begin{pmatrix}
        \widetilde{A} & C^T                           \\
        0             & -\frac{1}{\gamma}h^{\alpha} M
    \end{pmatrix}.
\end{equation}

We derive lower and upper bounds for the eigenvalues of the preconditioned matrix and then investigate the behavior of the smallest eigenvalue as the mesh is refined.

While eigenvalue information alone is generally insufficient to fully characterize the convergence of non-symmetric Krylov subspace methods, especially when the problem is far from normal~\cite{GMRES_eigs}, rapid convergence is commonly observed when the spectrum is real and positive and lies within a bounded interval away from zero. The spectral analysis below is performed for the ideal preconditioner~\eqref{eqn:candidate_prec}. However, we anticipate using a more practical approximation thereof in the numerical experiments.

From now on, we will assume $a_0(\cdot, \cdot)$ to be not only positive semidefinite, but also symmetric, and $C$ to have full row rank, consistently with the inf-sup condition~\eqref{eqn:inf-sup}.

Note that, in principle, the natural choice of $W$ in~\eqref{eqn:gen_aug_sys} should be consistent with the norm on $\Lambda$ used in the augmented formulation~\eqref{eqn:problem_general_discrete_augmented_Lambda}. However, since we move from the norm on $\Lambda$ to that on $H$ in~\eqref{eqn:problem_general_discrete_augmented}, a mesh-dependent scaling factor arises in the choice of this matrix. To emphasize this correspondence, we set $\widetilde{W}\coloneqq h^{\alpha} M$ and rewrite the matrix in~\eqref{eqn:aug_system_matrix} and the preconditioner in~\eqref{eqn:candidate_prec} accordingly.

\begin{theorem}\label{thm:lambda_bound}
    Assume that $\widetilde{\mathcal{A}}_\gamma$ and $\widetilde{\mathcal{P}}_{\gamma}$ are defined by the matrices in Equations~\eqref{eqn:aug_system_matrix} and~\eqref{eqn:candidate_prec}, respectively. The eigenvalues of
    the preconditioned matrix $\widetilde{\mathcal{P}}_{\gamma}^{-1} \widetilde{\mathcal{A}}_\gamma$ are all real and positive. %
    It holds $$\operatorname{Spec}\Bigl(\widetilde{\mathcal{P}}_{\gamma}^{-1} \widetilde{\mathcal{A}}_\gamma\Bigr) \subseteq \left [ \eta, 1 \right ], $$
    where
    $$ \eta \coloneqq \min_{\mathbf{0} \neq \mathbf{v} \in \mathcal{S}}  \frac{\gamma \mathbf{v}^T C^T \Winv C \mathbf{v}}{\mathbf{v}^T \hat{A} \mathbf{v} + \gamma \mathbf{v}^T C^T \Winv C \mathbf{v}},
    $$
    \noindent
    and $\nu=1$ is an eigenvalue of algebraic multiplicity at least $n$. Here $\mathcal{S}$ denotes the set
    \begin{equation}\label{eqn:S_set}
        \mathcal{S}\coloneqq\{\mathbf v \in \mathbb{R}^{n} \mid \mathbf{u}^T \hat{A} \mathbf{v} =0\>\> \forall \> \mathbf{u} \in \ker(C)\}.
    \end{equation}
\end{theorem}

\vspace{0.1cm}
\begin{proof}
    Since both the original system matrix $\widetilde{\mathcal{A}}_\gamma$ and the preconditioner are nonsingular, the preconditioned matrix is also nonsingular, and therefore $\nu=0$ cannot be an eigenvalue.
    Let $\nu\neq 0$ be an arbitrary eigenvalue of the preconditioned matrix, with a corresponding
    eigenvector $(\mathbf{x;y})$.
    The generalized eigenvalue problem can be stated as:
    \begin{equation*}
        \begin{pmatrix}
            \hat{A} + \gamma C^T \Winv C & C^T \\
            C                            & 0
        \end{pmatrix}
        \begin{pmatrix}
            \mathbf{x} \\
            \mathbf{y}
        \end{pmatrix}= \nu
        \begin{pmatrix}
            \hat{A} + \gamma C^T \Winv C & C^T                 \\
            0                            & -\frac{1}{\gamma}\W
        \end{pmatrix}
        \begin{pmatrix}
            \mathbf{x} \\
            \mathbf{y}
        \end{pmatrix},
    \end{equation*}
    which can be written explicitly as the following system of equations:
    \begin{align}
        (\hat{A} + \gamma C^T \Winv C) \mathbf{x} + C^T \mathbf{y} & = \nu \big((\hat{A} + \gamma C^T \Winv C) \mathbf{x} + C^T \mathbf{y}\big), \label{eqn:first_row_eigenproblem} \\
        C\mathbf{x}                                                & = -\frac{\nu}{\gamma}\W \mathbf{y}.\label{eqn:second_row_eigenproblem}
    \end{align}

    Notice that $\mathbf{x} \ne \mathbf{0}$; otherwise, the positive definiteness of $\W$ implies $\mathbf{y}=\mathbf{0}$, in
    contradiction with the fact that $(\mathbf{x},\mathbf{y})$ is an eigenvector. Therefore, we can assume $\mathbf{x} \ne \mathbf{0}$.

    It is evident from Equation~\eqref{eqn:first_row_eigenproblem} that $\nu=1$ is an eigenvalue of
    $\widetilde{\mathcal{P}}_{\gamma}^{-1} \widetilde{\mathcal{A}}_\gamma$, with
    associated eigenvector $(\mathbf{x};-\gamma \Winv C\mathbf{x})$.
    Moreover, when $\mathbf{x} \in \ker(C)$, we
    have $(\mathbf{x};\mathbf{0})$ as associated eigenvector.

    We now assume $\nu \ne 1$. From~\eqref{eqn:first_row_eigenproblem} we derive
    \begin{equation} \label{eqn:only_x}
        (\hat{A} + \gamma C^T \Winv C) \mathbf{x} + C^T \mathbf{y}  = 0,
    \end{equation}
    whereas from~\eqref{eqn:second_row_eigenproblem} we obtain $\mathbf{y}=-\frac{\gamma}{\nu}\Winv C \mathbf{x}$, which plugged into the previous equation
    gives
    \begin{equation}\label{eqn:first_row_mod}
        (\hat{A} + \gamma C^T \Winv C)\mathbf{x} - \frac{\gamma}{\nu} C^T \Winv C \mathbf{x} = 0.
    \end{equation}
    Multiplying both sides of~\eqref{eqn:first_row_mod} by $\nu  \mathbf{x}^\ast$, we get
    \begin{equation}\label{ref:eq_for_lambda}
        \nu \mathbf{x}^\ast(\hat{A} + \gamma C^T \Winv C) \mathbf{x} - \gamma \mathbf{x}^\ast C^T \Winv C \mathbf{x} = 0.
    \end{equation}
    Therefore,
    \begin{equation}\label{eqn:lambda}
        \nu  =  \frac{\gamma \mathbf{x}^\ast C^T \Winv C \mathbf{x}}{\mathbf{x}^\ast \hat{A}  \mathbf{x} + \gamma {\mathbf{x}^\ast C^T \Winv C \mathbf{x}}}
        =\frac{\gamma \mathbf{x}^\ast C^T \Winv C \mathbf{x}}{\mathbf{x}^\ast A  \mathbf{x} + \gamma h^{-\alpha}\mathbf{x}^{\ast}E\mathbf{x} + \gamma {\mathbf{x}^\ast C^T \Winv C \mathbf{x}}}.
    \end{equation}

    Since $C$ has full rank, a necessary and sufficient condition for the saddle point matrix $\widetilde{\mathcal{A}}_\gamma$ to be nonsingular is $\ker(\widetilde A) \cap \ker(C) = \{\mathbf{0}\}$, which implies $\ker(\hat A) \cap \ker(C) = \{\mathbf{0}\}$. Therefore, the augmented matrix $\hat{A} + \gamma C^T \Winv C$ is still symmetric positive definite even if $A$ is SPSD.
    The positive definiteness of $\hat A + \gamma C^T \Winv C$ and the positive semidefiniteness of $C^T \Winv C$ imply that all the eigenvalues are positive and real\footnote{Hence, the corresponding eigenvector can also be chosen to be real: for this reason $\mathbf{x}^\ast$ will be replaced by $\mathbf{x}^T$.}.
    Moreover, note that in this case $\mathbf{x}$ cannot belong to $\ker(\hat{A})$, since $\nu \neq 1$, nor to $\ker(C)$ since $\nu \neq 0$.

    From~\eqref{eqn:lambda}, we also deduce that $\nu < 1$, and that all eigenvalues satisfying~\eqref{eqn:lambda} cluster towards $1$ as
    $\gamma \rightarrow +\infty$ when the augmentation does not contain the term $\gamma h^{-\alpha}E$. In view of Lemmas~\ref{lemma:algebraic_term} and~\ref{lemma:stab_term}, this is precisely the case for algebraic augmentation, for which no stabilizing term is present, consistently with the theory of classical AL-based preconditioners~\cite{ALprec}.

    To conclude, we know that all nonunit eigenvalues of $\widetilde{\mathcal{P}}_{\gamma}^{-1} \widetilde{\mathcal{A}}_\gamma$ satisfy~\eqref{ref:eq_for_lambda}, hence, by the theory of generalized Rayleigh quotients, the smallest positive eigenvalue is characterized by:
    $$ \eta \coloneqq \min_{\mathbf{0} \neq \mathbf{v} \in \mathcal{S}}  \frac{\gamma \mathbf{v}^T C^T \Winv C \mathbf{v}}{\mathbf{v}^T \hat{A} \mathbf{v} + \gamma \mathbf{v}^T C^T \Winv C \mathbf{v}},
    $$ which concludes the proof.
\end{proof}
\begin{remark}
    There are particular cases in which the matrix $\widetilde{\mathcal{A}}_\gamma$ is actually singular, for example, in the presence of enclosed flows in the Stokes problem, where $a_0(\cdot, \cdot)$ is positive definite but the divergence matrix $C$ is rank deficient by one. However, the presence of a zero eigenvalue does not necessarily hinder the convergence of preconditioned GMRES, provided that the system is consistent~\cite{modALprec,ElmanSilvesterWathen}.
\end{remark}

\subsection{Mesh-independence of the lower bound}
The behavior of $\eta$ upon mesh refinement relies on the algebraic formulation of the inf-sup condition~\eqref{eqn:inf-sup}.

Whenever $\Lambda \subseteq H$, we assume that the multiplier space satisfies an inverse estimate of the form
\begin{equation}\label{eqn:inequality_multiplier}
    \|\mu_h\|_{\Lambda'} \geq c h^{s} \|\mu_h\|_{H}
    \qquad \forall \mu_h \in \Lambda_h,
\end{equation}
where $h$ is the mesh size of the grid over which the multiplier space is defined.

Practical examples of such estimates will be discussed later in Remark~\ref{rmk:exponent}. Note that the case $\Lambda \equiv H$ can be interpreted as the limiting case corresponding to $s=0$, as no inverse estimate is needed.

Furthermore, let $\{\mathbf{v}_i\}_{i=1}^n$ and $\{\vb{\mu}_k\}_{k=1}^m$ be the expansion coefficients of $v_h$ and $\mu_h$ with respect to their bases, and
let $T$ denote the matrix corresponding to the inner product on $V$, so that $\|v_h\|_{V}^2 =  \mathbf{v}^T T \mathbf{v}$. After exploiting inequality~\eqref{eqn:inequality_multiplier}, we have the following algebraic interpretation of the discrete inf-sup stability condition~\eqref{eqn:inf-sup}:
\begin{equation}\label{eq:alg_inf-sup}
    \inf_{\vb{0}\ne \vb{\mu} } \sup_{\vb{0} \ne \mathbf{v}} \frac{\mathbf{v}^T C^T \vb{\mu}}{(\mathbf{v}^T T \mathbf{v})^{\frac{1}{2}}(\vb{\mu}^T h^{2s} M \vb{\mu})^{\frac{1}{2}}}\ge  \tilde\beta_h,
\end{equation}
where the realization of the dual norm has been replaced by the matrix $h^{2s} M$, which represents the mesh-dependent
norm $h^{2s}\|\cdot\|_H^2$. It is well known that the constant $\tilde \beta_h$ can be identified through a generalized eigenvalue problem (see, e.g.,~\cite{malkus, boffi2013mixed, ElmanSilvesterWathen} for details).
\begin{proposition}\label{prop:malkus}
    Let $\tilde \beta_h$ be the inf-sup constant in~\eqref{eq:alg_inf-sup}, and let $\nu_{\min}$ be the smallest positive generalized eigenvalue of
    $$C^T (h^{2s} M)^{-1} C \mathbf{v}= \nu T \mathbf{v}.$$
    Then, $$\tilde \beta_h^2 = {\nu_{\min}}= \min_{\mathbf{0} \neq \mathbf{v} \in \mathcal{R}}  \frac{ \mathbf{v}^T C^T (h^{2s} M)^{-1} C \mathbf{v}}{\mathbf{v}^T T \mathbf{v} },
    $$where $\mathcal{R}$ denotes the set $\mathcal{R}\coloneqq\{\mathbf v \in \mathbb{R}^{n} \mid \mathbf{u}^T T \mathbf{v} =0\>\> \forall \> \mathbf{u} \in \ker(C)\}$.

\end{proposition}

In the next theorem, we show that requiring the smallest eigenvalue $\eta$ to be independent of the mesh size leads to a relation between the exponent $\alpha$ in the augmented formulation~\eqref{eqn:problem_general_discrete_augmented} and the exponent $s$ appearing in inequality~\eqref{eqn:inequality_multiplier}. In addition, we quantify the influence of the stabilizing term when it is nonzero. To this aim, we define
\begin{equation}\label{eqn:ckappa}
    c_\kappa(\mathbf{v}) \coloneqq \frac{h^{-2s} \mathbf{v}^T E \mathbf{v}}{\mathbf{v}^ TT\mathbf{v}} = \frac{h^{-2s}\|\kappa(L v_h)\|_H^2}{\|v_h\|_V^2}.
\end{equation}

\begin{theorem}\label{thm:lower_bound}
    Let $\gamma$ be a positive constant, and assume that the inf-sup condition~\eqref{eqn:inf-sup} holds and that the multiplier space satisfies the inverse estimate~\eqref{eqn:inequality_multiplier}. Then, for $\alpha=2s$, the lower bound $\eta$ for the spectrum of $\widetilde{\mathcal{P}}_{\gamma}^{-1} \widetilde{\mathcal{A}}_\gamma$ in Theorem~\ref{thm:lambda_bound} satisfies
    \begin{equation}\label{eqn:lower_bound}
        \eta \coloneqq \min_{\mathbf{0} \neq \mathbf{v} \in \mathcal{S}}  \frac{\gamma \mathbf{v}^T C^T \Winv C \mathbf{v}}{\mathbf{v}^T \hat{A} \mathbf{v} + \gamma \mathbf{v}^T C^T \Winv C \mathbf{v}}
        \geq
        \min_{\mathbf{0} \neq \mathbf{v} \in \mathcal{S}}
        \frac{1}{1+\dfrac{c_a}{\gamma\,\tilde \beta_h^2 } + \dfrac{ c_{\kappa}(\mathbf{v}_\perp)}{\tilde \beta_h^2}
        }
    \end{equation}where $c_a$ denotes the continuity constant of $a_0(\cdot,\cdot)$ and $\mathbf{v}_\perp$ is the component of $\mathbf{v}$ that is T-orthogonal to $\ker(C)$. In particular, if $c_\kappa(\cdot)$ is uniformly bounded in $h$, then $\eta$ is bounded away from $0$ as $h \to 0$.
\end{theorem}
\begin{proof}
    We first observe that $\mathcal{S}\cap\ker(C)=\{\mathbf{0}\}$. Indeed, let $\mathbf{v}\in\mathcal{S}\cap\ker(C)$. Taking
    $\mathbf{u}=\mathbf{v}$ in the definition of $\mathcal{S}$ in~\eqref{eqn:S_set},
    we obtain $\mathbf{v}^T\hat{A}\mathbf{v}=0$.
    Since $\hat{A}$ is symmetric positive semidefinite, it follows that
    $\hat{A}\mathbf{v}=\mathbf{0}$,
    that is, $\mathbf{v}\in\ker(\hat{A})$. Hence, $\mathbf{v}\in\ker(\hat{A})\cap\ker(C)=\{\mathbf{0}\}$,
    which proves that $\mathcal{S}\cap\ker(C)=\{\mathbf{0}\}$.

    Now let $\mathbf{v}\in\mathcal{S}$ be nonzero. In order to use the results in Proposition~\ref{prop:malkus}, we decompose $\mathbf{v}$ in a $T$-orthogonal fashion as
    \begin{equation*}
        \mathbf{v}=\mathbf{v}_\perp+\mathbf{v}_0,\qquad
        \mathbf{v}_0\in\ker(C),\quad \mathbf{v}_\perp^T T\,\mathbf{v}_0=0,
    \end{equation*}
    so that $\mathbf{v}_\perp$ is $T$-orthogonal to $\ker(C)$ and
    $\mathbf{v}^T T\mathbf{v}=\mathbf{v}_\perp^T T\mathbf{v}_\perp+\mathbf{v}_0^T T\mathbf{v}_0$.
    Note that $\mathbf{v}_\perp\neq \mathbf{0}$, otherwise $\mathbf{v}=\mathbf{v}_0\in\ker(C)$.

    Since $\W$ is SPD and $\mathbf{v}_0\in\ker(C)$, the augmented term only sees the
    component $\mathbf{v}_\perp$:
    \begin{equation*}
        C^T \Winv C\,\mathbf{v}
        = C^T \Winv C\,\mathbf{v}_\perp.
    \end{equation*}
    As $\mathbf{v}_\perp$ is $T$-orthogonal to $\ker(C)$, it belongs to $\mathcal{R}$. Therefore,
    Proposition~\ref{prop:malkus} yields
    \begin{equation}\label{eqn:numerator_bound}
        \mathbf{v}_\perp^T C^T M^{-1} C\,\mathbf{v}_\perp
        \geq \tilde \beta_h^2 h^{2s}\,\mathbf{v}_\perp^T T\,\mathbf{v}_\perp .
    \end{equation}
    Moreover, the defining property of $\mathcal{S}$ gives
    $\mathbf{v}_0^T \hat{A}\,\mathbf{v}=0$, whence $\mathbf{v}_0^T \hat{A}\,\mathbf{v}_\perp=-\,\mathbf{v}_0^T \hat{A}\,\mathbf{v}_0$.
    Using the symmetry and positive semidefiniteness of $\hat A$ together with the continuity of
    $a_0(\cdot,\cdot)$, we obtain
    \begin{equation}\label{eqn:Aterm_bound}
        \mathbf{v}^T \hat{A}\,\mathbf{v}
        =\mathbf{v}_\perp^T \hat{A}\,\mathbf{v}_\perp-\mathbf{v}_0^T \hat{A}\,\mathbf{v}_0
        \leq \mathbf{v}_\perp^T \hat{A}\,\mathbf{v}_\perp
        \leq c_a\,\mathbf{v}_\perp^T T\,\mathbf{v}_\perp + \gamma h^{-\alpha} \mathbf{v}_\perp^T E\,\mathbf{v}_\perp,
    \end{equation}
    where $c_a$ is the continuity constant for $a_0(\cdot,\cdot)$. Using the definition of $\W$, rewriting $\eta$ in a more convenient form and using~\eqref{eqn:numerator_bound} and~\eqref{eqn:Aterm_bound} together with the definition of $c_\kappa(\mathbf{v})$ in~\eqref{eqn:ckappa} we obtain:
    \begin{align*}
        \eta & \coloneqq \min_{\mathbf{0} \neq \mathbf{v}\in\mathcal{S}}
        \frac{\gamma\,\mathbf{v}^T C^T h^{-\alpha} M^{-1} C\,\mathbf{v}}
        {\mathbf{v}^T \hat{A} \mathbf{v}+\gamma \mathbf{v}^T C^T h^{-\alpha} M^{-1} C\,\mathbf{v}}
        = \min_{\mathbf{0} \neq \mathbf{v}\in\mathcal{S}}
        \frac{1}{1+\dfrac{\mathbf{v}^T \hat{A}\,\mathbf{v}}
        {\gamma\,\mathbf{v}^T C^T h^{-\alpha} M^{-1} C\,\mathbf{v}}}     \\
             & \geq \min_{\mathbf{0} \neq \mathbf{v}\in\mathcal{S}}
        \frac{1}{1+\dfrac{c_a \mathbf{v}_\perp^T T\,\mathbf{v}_\perp + \gamma h^{-\alpha}\mathbf{v}_\perp^T E\,\mathbf{v}_\perp\,}
            {\gamma\,\tilde \beta_h^2 h^{(2s-\alpha)} \,\mathbf{v}_\perp^T T\,\mathbf{v}_\perp}}
        \underbrace{=}_{\alpha=2s} \min_{\mathbf{0} \neq \mathbf{v}\in\mathcal{S}} \frac{1}{1+\dfrac{c_a}{\gamma\,\tilde \beta_h^2 } + \dfrac{c_\kappa(\mathbf{v}_\perp)}{\tilde \beta_h^2}}.%
    \end{align*}

    For algebraic augmentation, the term involving $c_\kappa(\mathbf{v}_\perp)$ is absent. Hence, for $\alpha=2s$, the lower bound in~\eqref{eqn:lower_bound} is independent of $h$, since $\tilde\beta_h$ is the mesh-independent inf-sup constant. Moreover, in agreement with Theorem~\ref{thm:lambda_bound}, increasing the penalization parameter $\gamma$ drives $\eta$ closer to $1$.

    For continuous augmentation, on the other hand, additional problem-dependent assumptions on $c_\kappa(\cdot)$ are generally required to establish a mesh-independent lower bound for $\eta$.
\end{proof}

\begin{remark}\label{rmk:exponent}
    The previous analysis relates the choice for the exponent $\alpha$ in the augmented formulation~\eqref{eqn:problem_general_discrete_augmented} to the power in the inverse estimate~\eqref{eqn:inequality_multiplier} on the multiplier space.
    This estimate depends on the problem, and it has to be verified on a case-by-case basis. We briefly list here some examples, which will also be considered in the numerical experiments.

    When $\Lambda = H^{\frac12}(\Gamma)$ and $H=L^2(\Gamma)$, where $\Gamma$ denotes either the boundary of the domain or an immersed interface, the inverse estimate~\eqref{eqn:inequality_multiplier} holds
    with $s=\frac12$. Consequently, $\alpha=1$ and $\W=h_\Gamma M$, where $h_\Gamma$ denotes the mesh size of the codimension-one manifold $\Gamma$.
    This setting arises, for instance, in the weak imposition of Dirichlet boundary conditions and in fictitious domain formulations with \textit{boundary} Lagrange multipliers~\cite{Glowinski}.

    When $\Lambda = H^{1}(\Omega_2)$ and $H=L^2(\Omega_2)$, where $\Omega_2$ is a subdomain of $\Omega$, the inverse estimate~\eqref{eqn:inequality_multiplier} holds with $s=1$, yielding
    $\alpha=2$ and $\W=h_{2}^2M$,
    where $h_{2}$ denotes the mesh size associated with the codimension-zero immersed domain $\Omega_2$.

    If $\Lambda \equiv H$ (as in the case of Grad-Div stabilization where both spaces coincide with $L^2(\Omega)$), no inverse estimate is needed. Therefore,  $s=0=\alpha$, and $\W = M$, the standard pressure mass matrix. More generally, this rule determines the exponent for the applications in Section~\ref{sec:numerical_tests}.
\end{remark}

We now validate the preceding framework on a concrete yet simple example involving the weak imposition of essential boundary conditions by means of Lagrange multipliers.

\subsection{Example: boundary conditions through Lagrange multiplier}\label{subsec:weak_bcs}
We consider the weak imposition of essential boundary conditions through Lagrange multipliers on a convex domain $\Omega \subset \mathbb{R}^d$ ($d=2,3$), introduced in Babu\v{s}ka's seminal work~\cite{BabuskaLM}.

Let $\Gamma = \partial \Omega$ and $\Lambda=H^{\frac{1}{2}}(\Gamma)$. Let $g \in H^{\frac{1}{2}}(\Gamma)$ be the prescribed boundary datum and $f \in L^2(\Omega)$
the forcing term. The operator $L : V \to \Lambda$, with $V = H^1(\Omega)$, is the trace operator, namely $Lu = u_{|\Gamma}$. The saddle-point formulation in~\cite{BabuskaLM} seeks $(u,\lambda) \in V \times \Lambda'$ such that
\begin{equation}\label{eqn:Babuska_LM}
    \begin{cases}
        (\nabla u, \nabla v)_{\Omega} + (u,v)_\Omega + \langle v,\lambda \rangle_{\Lambda \times \Lambda'}  = (f,v)_{\Omega}  \qquad  \forall v \in V, \\
        \langle u,\mu \rangle_{\Lambda \times \Lambda'}                                      = \langle g, \mu \rangle_{\Lambda \times \Lambda'} \qquad \qquad \qquad \quad \qquad \quad \>\> \forall \mu \in \Lambda'.
    \end{cases}
\end{equation}

Let us introduce suitable finite-dimensional spaces $V_h \subset V$ and $\Lambda_h \subset H \subset \Lambda'$, where $H=L^2(\Gamma)$. This choice is admissible since, identifying $L^2(\Gamma)$ with its dual, we have the dense and continuous inclusion
$$H^{\frac{1}{2}}(\Gamma) \subseteq L^2(\Gamma)  \subseteq H^{-\frac{1}{2}}(\Gamma).$$

For this problem, Remark~\ref{rmk:exponent} suggests $s=1/2$, and hence $\alpha=1$ as power of the mesh-size\footnote{\samepage Dimensional analysis also motivates $\alpha=1$. Assigning the unit $[U]$ to $u_h$, $v_h$, and $g$ and a characteristic length scale $L$, the terms $a(u_h,v_h)$ and $(u_h-g,v_h)_{\Gamma}$ scale as $U^2 L^{d-2}$ and $U^2 L^{d-1}$, respectively. Matching their units requires an augmentation coefficient scaling as $L^{-1}$, motivating the factor $h_{\Gamma}^{-1}$.}.

Let $h_{\Gamma}$ and $h_{\Omega}$ denote the mesh sizes on $\Gamma$ and $\Omega$, respectively, and assume that both meshes are shape regular and quasi-uniform. Furthermore, since $\Lambda_h \subset L^2(\Gamma)$, the pairings in~\eqref{eqn:Babuska_LM} reduce to $L^2(\Gamma)$ scalar products at the discrete level, so that the discrete augmented formulation reads as follows:
find $(u_h,\lambda_h) \in V_h \times \Lambda_h$ such that
\begin{equation}\label{eqn:LM_discrete_augmented}
    \begin{cases}
        \begin{aligned}[t]
             & (\nabla u_h, \nabla v_h)_{\Omega} + (u_h,v_h)_{\Omega} + \gamma h_{\Gamma}^{-1}(u_h,v_h)_{\Gamma} + (v_h,\lambda_h)_{\Gamma} \\
             & \qquad = (f,v_h)_{\Omega} + \gamma h_{\Gamma}^{-1}(g,v_h)_{\Gamma} \qquad \forall v_h \in V_h,
        \end{aligned} \\
        (u_h,\mu_h)_{\Gamma} = (g, \mu_h)_{\Gamma} \qquad \forall \mu_h \in \Lambda_h.
    \end{cases}
\end{equation}
Upon discretization, a classical choice is to use linear conforming elements for $V_h$ and $\Lambda_h$ (either $\mathcal{Q}^1$ or $\mathcal{P}^1$ elements, depending on the mesh type).
In~\cite{BabuskaLM}, it is shown that the corresponding discrete formulation, with equal order spaces, is inf-sup stable if $h_{\Gamma}/h_{\Omega} \ge c > 0$, with $c$ chosen sufficiently large.

Note that the augmented term involves the inner product on $H=L^2(\Gamma)$. Albeit simple, this example shows an instance of a mixed-dimensional problem where $\Lambda \ne H$.

As explained in~\cite{boffi1997analysis}, since the norm in $H$ is weaker than
the one in $\Lambda$, we ought to check the coercivity assumption in~\eqref{eqn:coercivity_augmented_eq} for the well-posedness of the problem, which reads now
$$\exists \delta >0 \colon  (\nabla v, \nabla v)_{\Omega} + (v,v)_{\Omega} +  \xi\|v\|_{0,\Gamma}^2  \geq \delta \|v\|_V^2 \>\> \>\>\forall v \in V,$$with $\xi >0$. %
The condition is satisfied automatically with $\delta=1$.

We start by showing that the bound in~\eqref{eqn:decay_stabilizing_term} tends to zero as the mesh is refined, implying that the proposed continuous augmentation is consistent. To this end, we adapt the following result from the error analysis in~\cite{boffi1997analysis}, setting $\gamma=1$ for ease of exposition. Note that, since $\gamma$ is a mesh-independent constant, this choice does not affect the following results. For simplicity, we assume that the mesh sizes are comparable, i.e., $h_\Gamma\simeq h_\Omega$, and denote both by $h$.

\begin{theorem}[\cite{boffi1997analysis}]\label{thm:a_priori_boffi}
    Assume the hypotheses of Section~\ref{sec:AL_formulations}. Let $(u, \lambda)$ be the solution of the original problem~\eqref{eqn:problem_general} and $(u_h, \lambda_h)$ be the solution of the discrete
    augmented problem~\eqref{eqn:problem_general_discrete_augmented}. Moreover, assume that $\lambda \in H$, $\alpha\ge0$, and $\gamma=1$. Then, the
    following error estimates hold:
    \[
        \|u-u_h\|_V  \leq c \Bigl( (1+  h^{-\frac{\alpha}{2}}) \inf_{v_h \in V_h} \|u-v_h\|_V + h^{\frac{\alpha}{2}} \inf_{\mu_h \in \Lambda_h} \|\lambda - \mu_h\|_H \Bigr),
    \]
    \[
        \|\lambda-\lambda_h\|_{\Lambda^{'}}  \leq c \Bigl( (1 +  h^{-\frac{\alpha}{2}} +  h^{-\alpha}) \inf_{v_h \in V_h} \|u-v_h\|_V +  \inf_{\mu_h \in \Lambda_h} \|\lambda - \mu_h\|_H \Bigr).
    \]
\end{theorem}

As observed in~\cite{boffi1997analysis}, the previous a priori estimates are not optimal, since $\|\lambda-\lambda_h\|_{\Lambda}$ is controlled by the $H$-norm of the interpolation error. In this specific case, such bounds are not sharp. Nevertheless, we use them to characterize the behavior of the stabilizing term. To satisfy the hypothesis $\lambda \in H$, we assume %
additional regularity $g\in H^{\frac{3}{2}}(\Gamma)$, so that by elliptic regularity $u\in H^2(\Omega)$. In particular, $\lambda=-\frac{\partial u}{\partial\boldsymbol{n}}\in H^{\frac{1}{2}}(\Gamma)$~\cite{BabuskaLM}, and hence $\lambda\in H$.

The following approximation properties hold for the chosen finite element spaces:
\[
    \inf_{v_h \in V_h} \|u-v_h\|_V  = \mathcal{O}(h),%
    \qquad
    \inf_{\mu_h \in \Lambda_h} \|\lambda - \mu_h\|_H = \mathcal{O}(h^{\frac{1}{2}}). %
\]
Using Theorem~\ref{thm:a_priori_boffi} with $\alpha=1$ and standard conforming interpolation estimates, yields the following a-priori error estimate for
the augmented formulation:
\[
    \|u-u_h\|_V  \leq c \Bigl( \bigl( h + h^{\frac{1}{2}}\bigr) \vert u\vert_{2,\Omega} + h\vert\lambda\vert_{\frac{1}{2},\Gamma}\Bigr).
\]
Moreover, for $v \in V$, it holds
\begin{equation}\label{eqn:bound_steinb}
    \|\kappa (L v)\|_H = \|\kappa (v|_{\Gamma})\|_H = \|v|_{\Gamma} - \pi_H (v|_{\Gamma})\|_H \leq c h^{\frac{1}{2}} |v|_{\frac{1}{2},\Gamma} \leq c h^{\frac{1}{2}} \|v\|_{1,\Omega},
\end{equation}where we have used the approximation property of the $L^2$-projection $\pi_H$ (see, e.g.,~\cite[Chapter 10]{Steinbach2007NumericalAM}) and the trace theorem. Exploiting twice the bound~\eqref{eqn:bound_steinb} in Lemma~\ref{lemma:stab_term} yields
$$|r(v_h)| \leq \gamma h^{-1}  \|\kappa (L (u_h- u))\|_H \|\kappa(L v_h)\|_H \leq \gamma  h^{-1} c^2 h \|u - u_h\|_V  \|v_h\|_V \rightarrow 0, \quad h \to 0.$$

We now turn our attention to the computational advantages of continuous over algebraic augmentation in terms of matrix sparsity and assembly cost. We provide numerical evidence of this matter through the following three-dimensional example.
\begin{example}[Computational gains]
    We set $\Omega \!=\! \mathcal{B}_{1}(0,0,0)$, the three-dimensional unit ball centered at the origin, and report in Figure~\ref{fig:sparsity_babuska_LM} the sparsity patterns and the number of nonzero entries (nnz) of the matrix $A$ associated with the original bilinear form $a_0(\cdot,\cdot)$, the matrix $\widetilde{A}$ associated with the augmented bilinear form, and the algebraically augmented $(1,1)$-block, $A_{\text{alg}} = A + C^T M^{-1} C$. The global size of the $(1,1)$-block is $3,817$. The original matrix $A$ and the matrix $\tilde{A}$ share the same sparsity pattern, whereas $A_{\text{alg}}$ is clearly
    much denser, by a factor larger than $16$. Although the augmentation considered
    here involves only a surface integral, the algebraic augmentation still leads to a dramatic increase in the number of nonzero entries.

    To further highlight the benefits of our approach, we compare the wall-clock times of the continuous and algebraic augmentations (the latter computed using the cheaper variant $C^T\operatorname{diag}(M)^{-1}C$) over a sequence of mesh refinements in Figure~\ref{fig:costs_comparison}.

    \begin{figure}[tbp]
        \centering
        \captionsetup[subfigure]{font=footnotesize,margin=0pt}
        \begin{subfigure}[b]{0.28\textwidth}
            \centering
            \includegraphics[width=0.9375\textwidth]{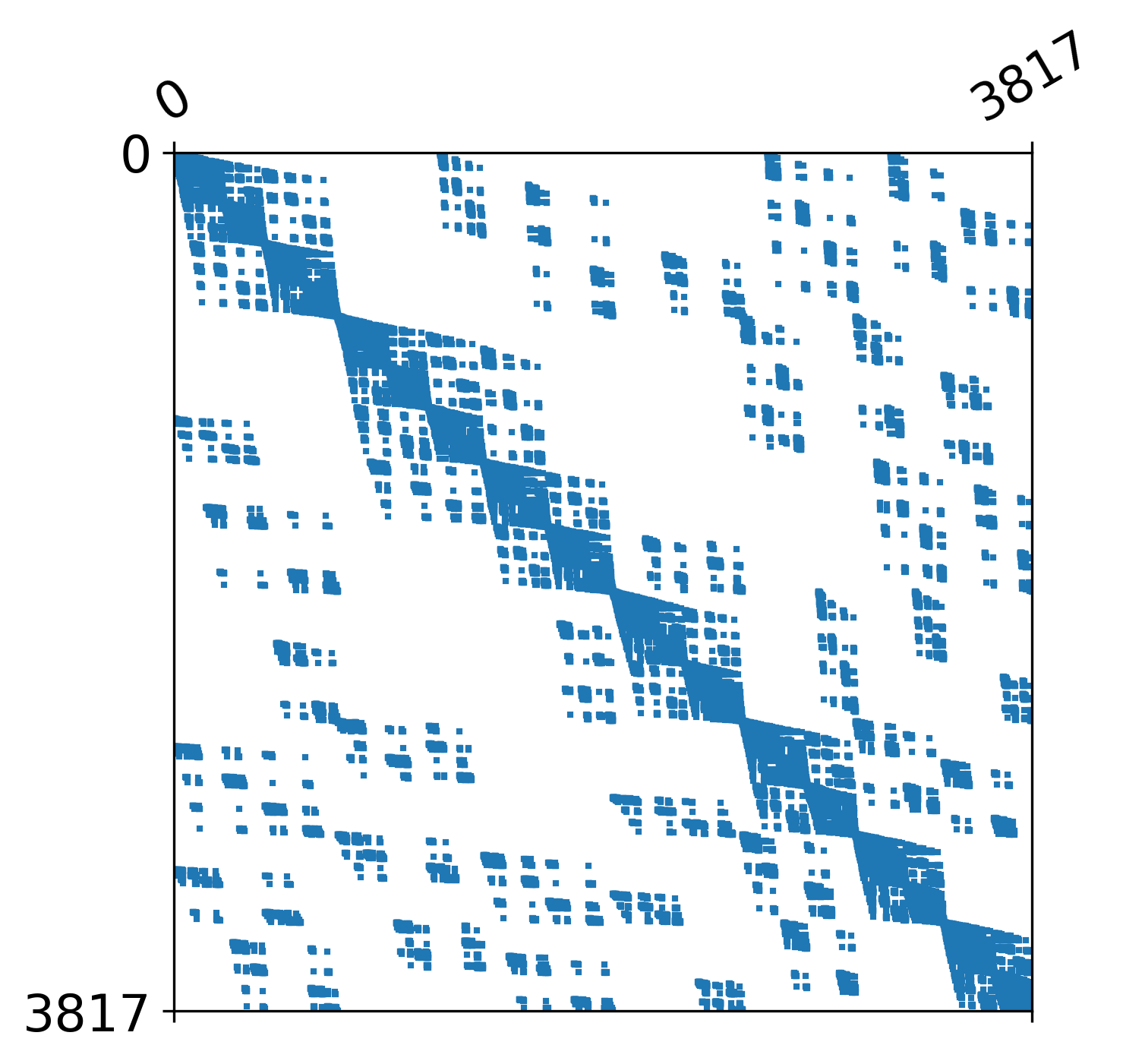}
            \caption{$A$ (nnz = 98,617)}
            \label{fig:sparsity:model1}
        \end{subfigure}
        \hfill
        \begin{subfigure}[b]{0.28\textwidth}
            \centering
            \includegraphics[width=0.9375\textwidth]{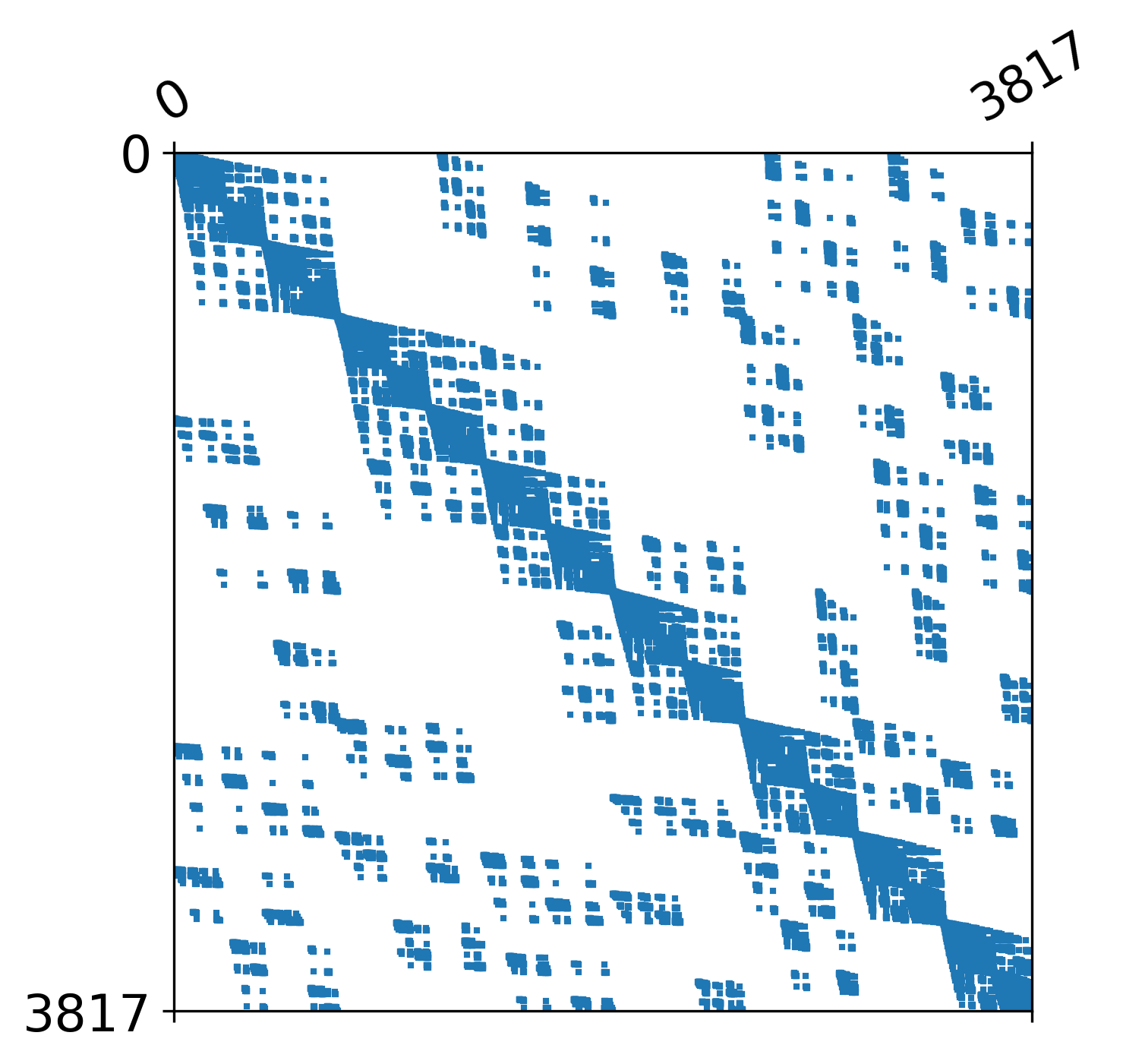}
            \caption{$\tilde{A}$ (nnz = 98,617)}
            \label{fig:sparsity:model2}
        \end{subfigure}
        \hfill
        \begin{subfigure}[b]{0.28\textwidth}
            \centering
            \includegraphics[width=0.9375\textwidth]{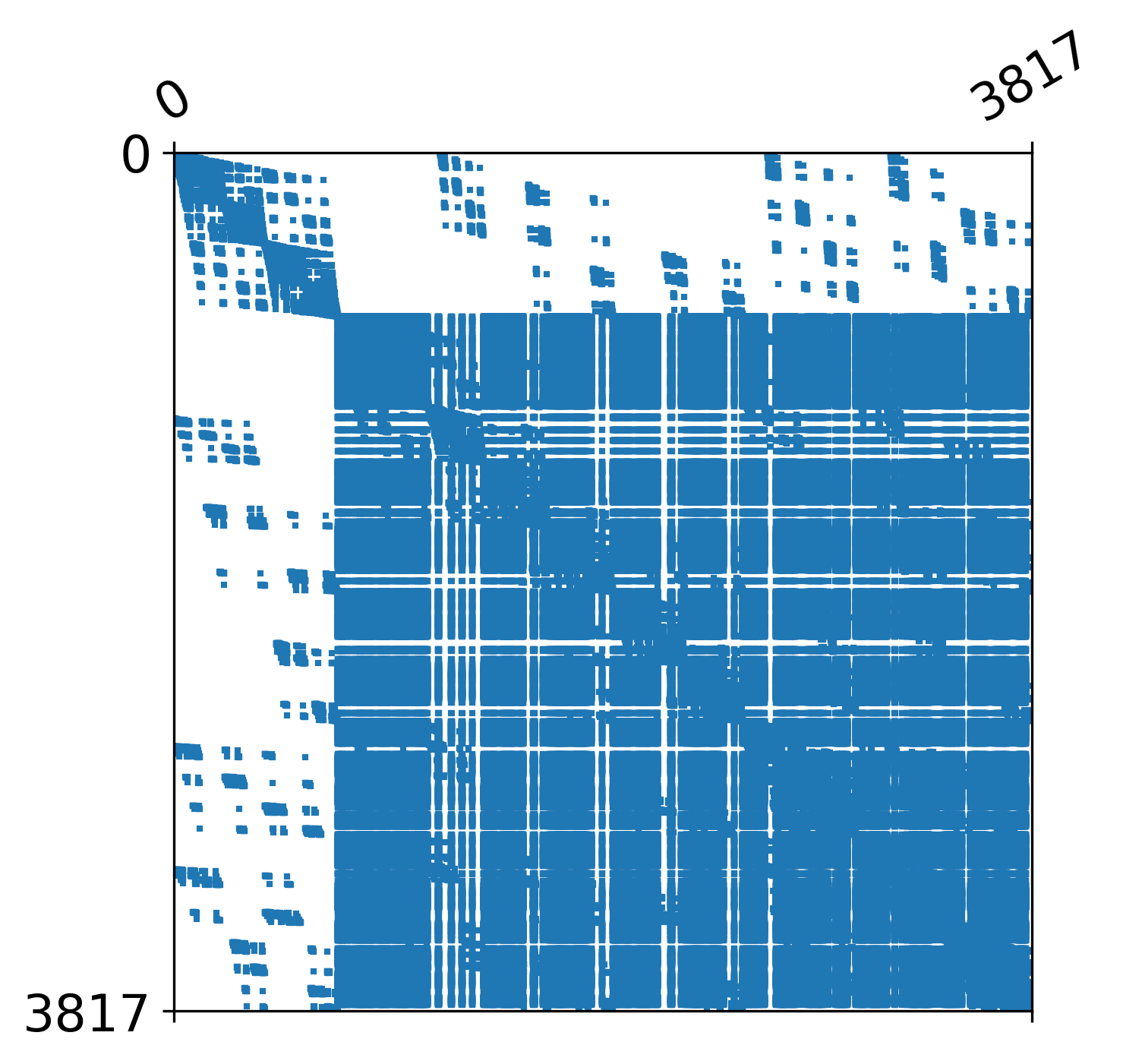}
            \caption{$A_{\text{alg}}$ (nnz=1,660,985)}
            \label{fig:sparsity:model3}
        \end{subfigure}
        \caption{Sparsity patterns of the original matrix $A$ and of the augmented matrices $\tilde{A}$, and $A_{\text{alg}}$ for the
            weak imposition of Dirichlet boundary conditions à la Babu\v{s}ka. $\mathcal{Q}^1$--$\mathcal{Q}^1$ elements are employed.}
        \label{fig:sparsity_babuska_LM}
    \end{figure}

    \begin{figure}[tbp]
        \centering
        \includegraphics{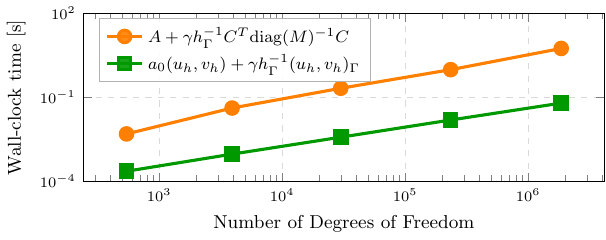}
        \caption{Comparison of wall-clock time (in seconds) against number of dof for the algebraic and continuous
            augmentations. Here $a_0(u_h, v_h) = (\nabla u_h, \nabla v_h)_{\Omega} + (u_h, v_h)_{\Omega}$, and $A$ is its
            associated matrix representation.}
        \label{fig:costs_comparison}
    \end{figure}
    The continuous augmentation is substantially sparser, and on the finest mesh is $90\times$ faster to assemble. Motivated by these examples, we will not compute any triple matrix product in the numerical
    tests of Section~\ref{sec:numerical_tests}, but rather we will directly assemble the augmented bilinear form.

    Finally, we establish the mesh-independence of $\eta$ in Theorem~\ref{thm:lower_bound} for the
    present problem, namely, we show that the term $c_\kappa(\mathbf{v})$ remains bounded as the mesh is refined. Indeed, using estimate~\eqref{eqn:bound_steinb} in the definition~\eqref{eqn:ckappa} yields $$ \frac{h^{-1}\|\kappa(L(v_h))\|_H^2 }{\|v_h\|_V^2} \leq c^2.$$
    \begin{figure}[!htbp]
        \centering
        \includegraphics[width=0.9\textwidth]{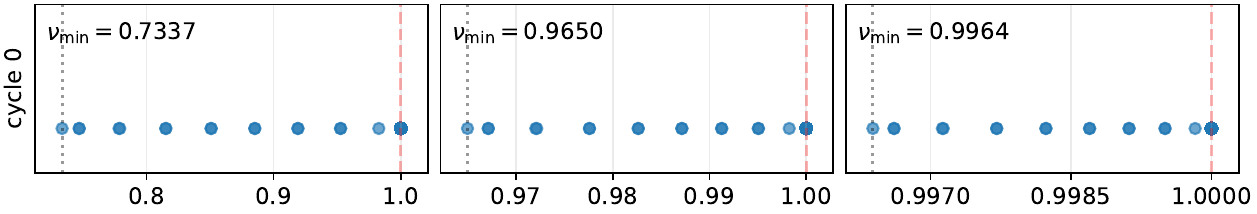}\vspace{0.5ex}
        \includegraphics[width=0.9\textwidth]{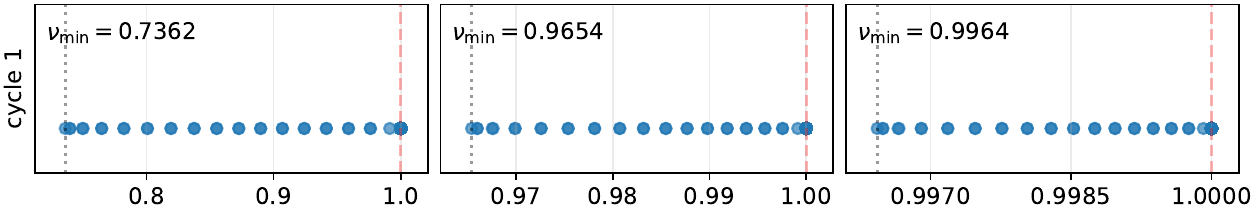}\vspace{0.5ex}
        \includegraphics[width=0.9\textwidth]{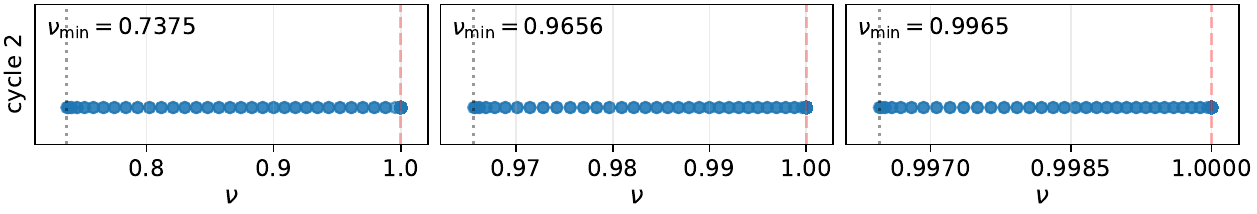}
        \caption{Spectrum of $\widetilde{\mathcal{P}}_{\gamma}^{-1} \widetilde{\mathcal{A}}_{\gamma}$ for the weak imposition of essential boundary conditions. Each row corresponds to a refinement cycle in which both meshes are uniformly refined. Along rows: influence of the AL parameter $\gamma \in \{10,100,1000\}$ (increasing from left to right). Along columns: mesh-independence of lower bound $\eta$. Cycle 0 has sizes $n=337$ and $m=16$, cycle 1 has $n=1,313$ and $m=32$, and cycle 2 has $n=5,185$ and $m=64$.}
        \label{fig:babuska_spectrum}
    \end{figure}

    We also validate the mesh-independence of $\eta$ numerically.
    To avoid solving large
    generalized eigenvalue problems, we consider the two-dimensional analogue of the previous example, setting $\Omega=\mathcal{B}_{1}(0,0)$.
    We perform a series of uniform refinement cycles to investigate the
    behavior of $\eta$ as the mesh resolution increases. Furthermore, we also
    verify that increasing the augmentation parameter using values $\gamma \in \{10,100,1000\}$ produces a tightly clustered spectrum near $1$, as predicted at the end of Theorem~\ref{thm:lower_bound}. We compute the eigenvalues of the preconditioned matrices using the Python \texttt{scipy.linalg} package. The plots in Figure~\ref{fig:babuska_spectrum} validate the developed theory: all eigenvalues
    are real and positive, with the smallest eigenvalue indicated in each plot. Fixing each column, we observe that the smallest eigenvalue remains uniformly bounded upon mesh refinement. Fixing each row, instead, we observe an
    increasingly clustered spectrum near $1$ as the constant $\gamma$ increases.
\end{example}

\subsubsection{Interplay between \texorpdfstring{$h^{-\alpha}$}{h\textasciicircum alpha}, convergence rates, and iteration counts}\label{subsubsec:weak_bc_alpha_power}

In~\cite{boffi1997analysis}, the power of $h$ is tuned to improve a-priori estimates, whereas in the present work $\alpha=1$ is dictated by preconditioning requirements. Since this choice is not driven by accuracy considerations, we have to verify that it does not deteriorate the convergence of the discrete solution.

To this end, we consider again Problem~\eqref{eqn:LM_discrete_augmented}, with $\Omega=(0,1)^2$, manufactured solution $u=\sin(\pi (x+0.3)) \sin(\pi (y+0.3))$, and $\mathcal{Q}^1$--$\mathcal{Q}^1$ elements. This pair is balanced, and the non-augmented scheme already achieves the best-approximation
rate~\cite{BabuskaLM}, so no accuracy can be gained from the augmentation. The purpose of the experiment is therefore to confirm that none is lost with $\alpha=1$, and to show that, in contrast, the performance of the iterative solver is strongly dependent on the choice of $\alpha$. We use FGMRES with an inexact variant of the ideal preconditioner to account for its variable application resulting from the iterative solution of the diagonal blocks. The absolute stopping tolerance is set to $10^{-12}$. The augmented block uses conjugate gradient (CG) with one V-cycle of algebraic multigrid and tolerance $10^{-2}$, while the mass matrix uses CG with point-Jacobi preconditioning. Details are provided in Section~\ref{sec:numerical_tests}. We set the augmentation parameter to $\gamma=10$.

Errors, convergence rates, and iteration counts are reported in Table~\ref{tab:weak_bc_alpha_comparison}, where $\alpha=1$ is compared with $\alpha=\frac{1}{2}$.

For both powers, the optimal convergence rates are confirmed (two in $L^2$ and one in semi-$H^1$), and the errors are essentially unaffected by $\alpha$. The iteration counts, instead, remain bounded only for $\W=h_\Gamma M$ (cf. also Table~\ref{tab:weak_bc_2d_convergence} in Section~\ref{sec:numerical_tests}).

\begin{table}[htbp]
    \centering
    \sisetup{group-separator={,}, group-minimum-digits=4}
    {\small
        \renewcommand{\arraystretch}{1.05}
        \resizebox{\linewidth}{!}{%
            \begin{tabular}{
                    S[table-format=7]
                    S[table-format=5]
                    c c c S[table-format=2]
                    c c c S[table-format=2]
                }
                \toprule
                        &
                        &
                \multicolumn{4}{c}{$\W = h_\Gamma^{\frac{1}{2}}M$}
                        &
                \multicolumn{4}{c}{$\W = h_\Gamma M$}
                \\
                \cmidrule(lr){3-6} \cmidrule(lr){7-10}
                {$\mathrm{dof}_u$}
                        & {$\mathrm{dof}_\lambda$}
                        & {$\|u-u_h\|_{0,\Omega}$}
                        & {$|u-u_h|_{1,\Omega}$}
                        & {$\|\lambda-\lambda_h\|_{0,\Gamma}$}
                        & {FGMRES its.}
                        & {$\|u-u_h\|_{0,\Omega}$}
                        & {$|u-u_h|_{1,\Omega}$}
                        & {$\|\lambda-\lambda_h\|_{0,\Gamma}$}
                        & {FGMRES its.}
                \\
                \midrule
                81      & 16                                   & $7.89 \times 10^{-3}$         & $2.60 \times 10^{-1}$         & $8.47 \times 10^{-1}$         & 10 & $7.66 \times 10^{-3}$         & $2.56 \times 10^{-1}$         & $8.47 \times 10^{-1}$         & 10 \\
                289     & 32                                   & $2.01 \times 10^{-3}$\,(1.97) & $1.31 \times 10^{-1}$\,(0.99) & $5.75 \times 10^{-1}$\,(0.56) & 14 & $1.92 \times 10^{-3}$\,(2.00) & $1.28 \times 10^{-1}$\,(1.00) & $5.75 \times 10^{-1}$\,(0.56) & 12 \\
                1089    & 64                                   & $5.12 \times 10^{-4}$\,(1.97) & $6.62 \times 10^{-2}$\,(0.99) & $4.02 \times 10^{-1}$\,(0.51) & 18 & $4.81 \times 10^{-4}$\,(2.00) & $6.39 \times 10^{-2}$\,(1.00) & $4.03 \times 10^{-1}$\,(0.51) & 12 \\
                4225    & 128                                  & $1.30 \times 10^{-4}$\,(1.97) & $3.34 \times 10^{-2}$\,(0.99) & $2.84 \times 10^{-1}$\,(0.50) & 21 & $1.20 \times 10^{-4}$\,(2.00) & $3.20 \times 10^{-2}$\,(1.00) & $2.84 \times 10^{-1}$\,(0.50) & 12 \\
                16641   & 256                                  & $3.31 \times 10^{-5}$\,(1.98) & $1.68 \times 10^{-2}$\,(0.99) & $2.00 \times 10^{-1}$\,(0.50) & 24 & $3.01 \times 10^{-5}$\,(2.00) & $1.60 \times 10^{-2}$\,(1.00) & $2.01 \times 10^{-1}$\,(0.50) & 12 \\
                66049   & 512                                  & $8.38 \times 10^{-6}$\,(1.98) & $8.48 \times 10^{-3}$\,(0.99) & $1.42 \times 10^{-1}$\,(0.50) & 27 & $7.51 \times 10^{-6}$\,(2.00) & $7.99 \times 10^{-3}$\,(1.00) & $1.42 \times 10^{-1}$\,(0.50) & 11 \\
                263169  & 1024                                 & $2.11 \times 10^{-6}$\,(1.98) & $4.27 \times 10^{-3}$\,(0.99) & $1.00 \times 10^{-1}$\,(0.50) & 30 & $1.88 \times 10^{-6}$\,(2.00) & $3.99 \times 10^{-3}$\,(1.00) & $1.00 \times 10^{-1}$\,(0.50) & 11 \\
                1050625 & 2048                                 & $5.34 \times 10^{-7}$\,(1.99) & $2.14 \times 10^{-3}$\,(0.99) & $7.10 \times 10^{-2}$\,(0.50) & 34 & $4.70 \times 10^{-7}$\,(2.00) & $2.00 \times 10^{-3}$\,(1.00) & $7.10 \times 10^{-2}$\,(0.50) & 11 \\
                \bottomrule
            \end{tabular}}}
    \caption{Convergence study for the weak imposition of boundary conditions on $\Omega=(0,1)^2$, $u=\sin(\pi (x+0.3)) \sin(\pi (y+0.3))$,
    $\mathcal{Q}^1$--$\mathcal{Q}^1$ elements, $\gamma=10$. Comparison between $\W=h_{\Gamma}^{\frac{1}{2}} M$ and the exponent, $\alpha=1$, predicted by Remark~\ref{rmk:exponent}, implying $\W=h_\Gamma M$. Convergence rates are
    reported in parentheses.}
    \label{tab:weak_bc_alpha_comparison}
\end{table}

\section{Numerical experiments}\label{sec:numerical_tests}

We evaluate the performance of the proposed preconditioners for linear systems stemming from mixed finite element discretizations of mixed-dimensional and fictitious domain formulations with Lagrange multipliers. Our previous works~\cite{BENZI2026118522,benzi2026augmentedlagrangianpreconditionersfictitious} demonstrated the robustness of algebraic AL approaches for fictitious domain problems, showing superior performance compared with other methods. Here, we focus on improving and generalizing the framework developed therein, while referring to our previous works for comparisons with alternative methods.

Computations have been carried out on a laptop equipped with an AMD Ryzen 9 PRO 7940HS processor and 64 GB of RAM. Experiments
have been performed using the \textsc{C++} finite element library \textsc{deal.II}~\cite{dealIIdesign}.

Our implementation allows exploiting several features, such as adaptive mesh refinement, efficient handling of non-matching meshes in parallel setting, as well as state-of-the-art algebraic multigrid preconditioners available in well-established linear algebra libraries. When we use conforming elements on quadrilaterals or hexahedra, local refinement induces hanging nodes. Continuity of the discrete solution is enforced by imposing constraints on the degrees of freedom (dof). The nodal unknowns associated with such nodes are eliminated by restricting them to an
interpolation of the nodal values of the finite element solution on the corresponding parent element~\cite{dealIIdesign}.

Table~\ref{tab:problems_preconditioners} summarizes the problems, preconditioners, and stable finite element pairs, including three-dimensional elasticity with different Lamé constants in each subdomain.

\begin{table}[htbp]
    \centering
    \setlength{\tabcolsep}{6pt}
    \renewcommand{\arraystretch}{1.20}
    \begin{tabular*}{\linewidth}{@{\extracolsep{\fill}} l c c @{}}
        \toprule
        \textbf{Problem} & \textbf{Preconditioner} & \textbf{Finite element space} \\
        \midrule
        Babu\v{s}ka LM   & $\widetilde{\mathcal{P}}_{\gamma}^{\text{Babu\v{s}ka}}$~\eqref{eqn:AL_prec_babuska} & $\mathcal{Q}^1$--$\mathcal{Q}^1$ (2D), $\mathcal{P}^1$--$\mathcal{P}^1$ (3D) \\
        FD-BLM           & $\widetilde{\mathcal{P}}_{\gamma}^{\mathrm{BLM}}$~\eqref{eqn:AL_prec_BLM} & $\mathcal{Q}^1$--$\mathcal{P}^1$ \\
        FD-DLM           & $\widetilde{\mathcal{P}}_{\gamma_1,\gamma_2}^{\mathrm{DLM}}$~\eqref{eqn:modified_AL_elliptic_interface_matrix_form} & $\mathcal{Q}^1$--$(\mathcal{Q}^1 + \mathfrak{B})$--$\mathcal{Q}^0$ \\
        FD-DLM (linear elasticity)           & $\widetilde{\mathcal{P}}_{\gamma_1,\gamma_2}^{\mathrm{DLM}}$~\eqref{eqn:modified_AL_elliptic_interface_matrix_form} & $\mathcal{Q}^1$--$(\mathcal{Q}^1 + \mathfrak{B})$--$\mathcal{Q}^0$ \\
        \bottomrule
    \end{tabular*}
    \caption{Problems considered in this section, corresponding AL preconditioners, and finite element spaces.}
    \label{tab:problems_preconditioners}
\end{table}

We use simplicial, quadrilateral, or hexahedral meshes according to the geometry. As usual in fictitious domain experiments, we use a structured background mesh $\mathcal{T}_h$ on $\Omega$.

The outer Krylov solver is set
to FGMRES~\cite{FGMRES} with a restart every 30 iterations. All iterations start from a zero initial guess and, unless otherwise specified, are stopped when the relative residual norm drops
below $10^{-9}$, or when the absolute residual is below $10^{-12}$. All preconditioners considered are applied inexactly. Inner solves for the augmented blocks use preconditioned conjugate gradient (PCG) with one V-cycle of an algebraic multigrid preconditioner from \textsc{Trilinos ML}~\cite{Trilinos} (using a degree two Chebyshev smoother, aggregation threshold set to $10^{-3}$, maximum size of the coarsest level set to $2000$, and a direct solver on the coarsest level)
with a loose inner tolerance of $10^{-2}$. For the mass matrix inversion we use PCG with point-Jacobi. In each one of the forthcoming tables, we report as last column the (average) number of PCG iterations to invert the (1,1)-block.

A classical alternative for the approximate inversion of augmented blocks in AL preconditioners is to use one or a few multigrid cycles with tailored smoothers based on patchwise relaxation~\cite{schoeberl_1999,pcpatch}.
Crucially, designing such smoothers strongly depends on the weak form and on the finite element spaces employed; we notice that, in the mixed-dimensional case, this becomes even more challenging due to the presence of two overlapping meshes. For these reasons, we do
not pursue this avenue here, and leave the development of tailored multigrid solvers to future work.

\subsection{Weak imposition of boundary conditions with Lagrange multipliers}
We consider the system associated with the weak imposition of boundary conditions on $\Gamma = \partial \Omega$ through Lagrange multipliers in Problem~\eqref{eqn:LM_discrete_augmented}. We use a stable pair with linear Lagrangian elements for both the solution variable and the Lagrange multiplier. The resulting
preconditioner is~\eqref{eqn:candidate_prec} with blocks \begin{equation}\label{eqn:AL_prec_babuska}
    \begin{aligned}
        \widetilde{A}_{ij} & = (\nabla \phi_i, \nabla \phi_j)_{\Omega} + (\phi_i, \phi_j)_{\Omega} + \gamma h_{\Gamma}^{-1} ( \phi_i, \phi_j )_{\Gamma}, \\
        M_{k\beta}         & = (\psi_k, \psi_\beta)_{\Gamma},                                                                                            \\
        C_{ki}             & = (\phi_i,\psi_k)_{\Gamma}.
    \end{aligned}
\end{equation}
Following Remark~\ref{rmk:exponent}, we have set $s=\frac{1}{2}$, which corresponds to $\W = h_{\Gamma}M$. We refer to this preconditioner as $\widetilde{P}_\gamma^{\text{Babu\v{s}ka}}$.
We test its performance by setting $\gamma=10$, with $\Omega=(0,1)^2$ and $g=x+y$ as a non-constant boundary datum.
The results are reported in Table~\ref{tab:weak_bc_2d_convergence}. They indicate excellent convergence properties of the AL preconditioner, which exhibits constant and very low iteration counts across mesh refinement, even with a mild value of $\gamma=10$.

\begin{figure}[tbp]
    \centering
    \begin{subfigure}[t]{0.175\linewidth}
        \centering
        \includegraphics[width=\linewidth]{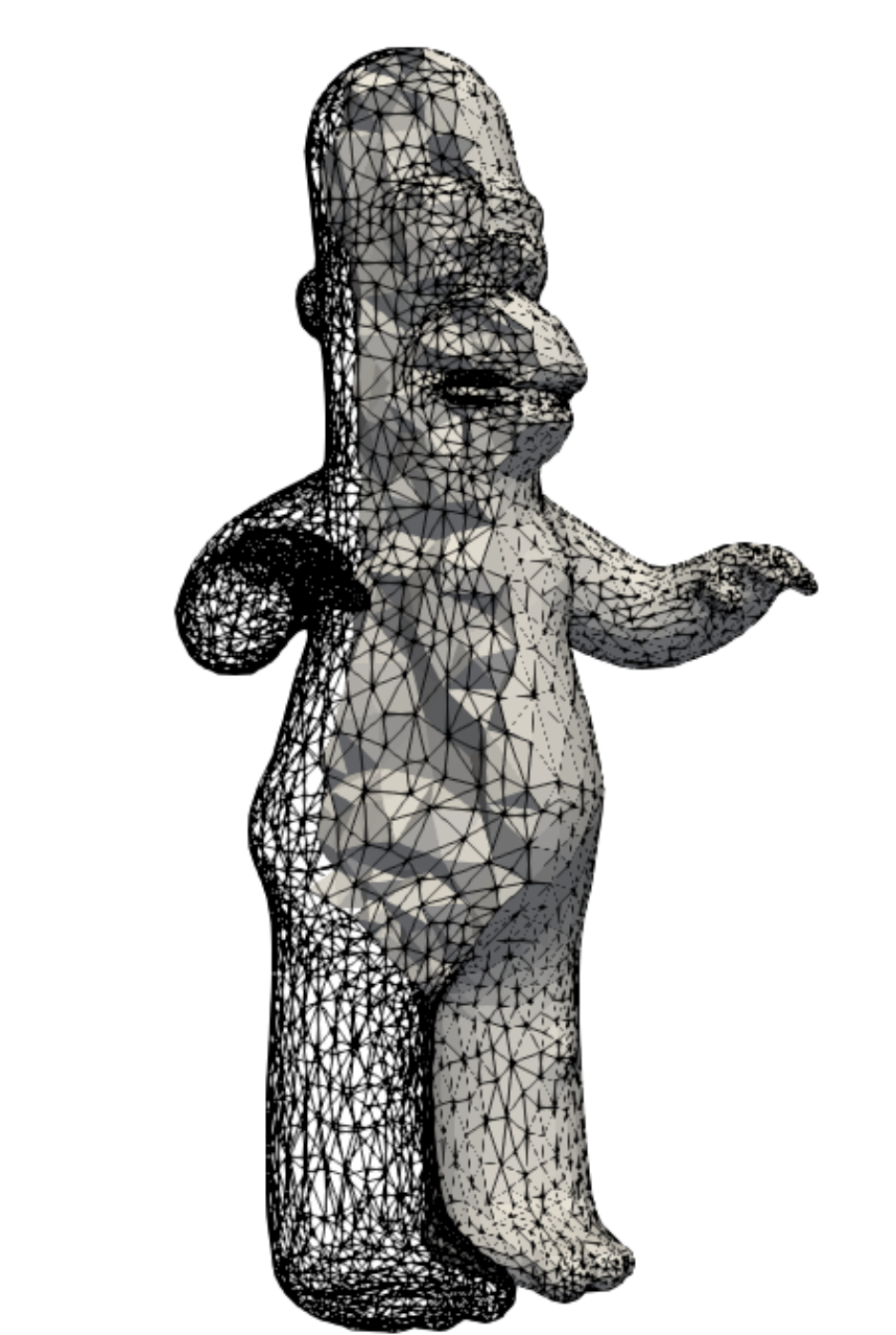}
        \caption{Mesh.}
        \label{fig:weak_bc_3d_geometry_shell}
    \end{subfigure}\hfill
    \begin{subfigure}[t]{0.175\linewidth}
        \centering
        \includegraphics[width=\linewidth]{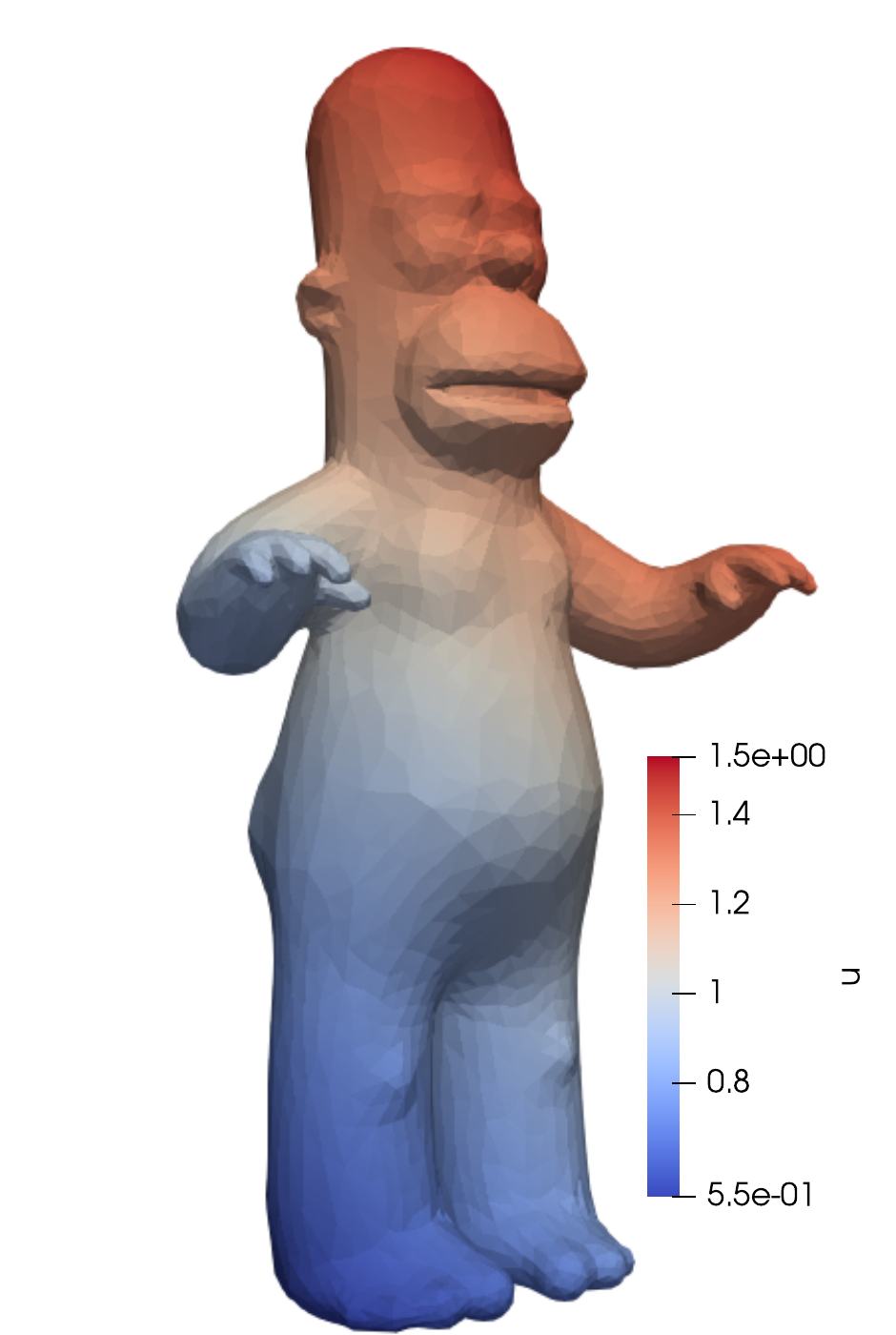}
        \caption{Solution.}
        \label{fig:weak_bc_3d_geometry_sol}
    \end{subfigure}\hfill
    \begin{subfigure}[t]{0.275\linewidth}
        \centering
        \includegraphics[width=\linewidth]{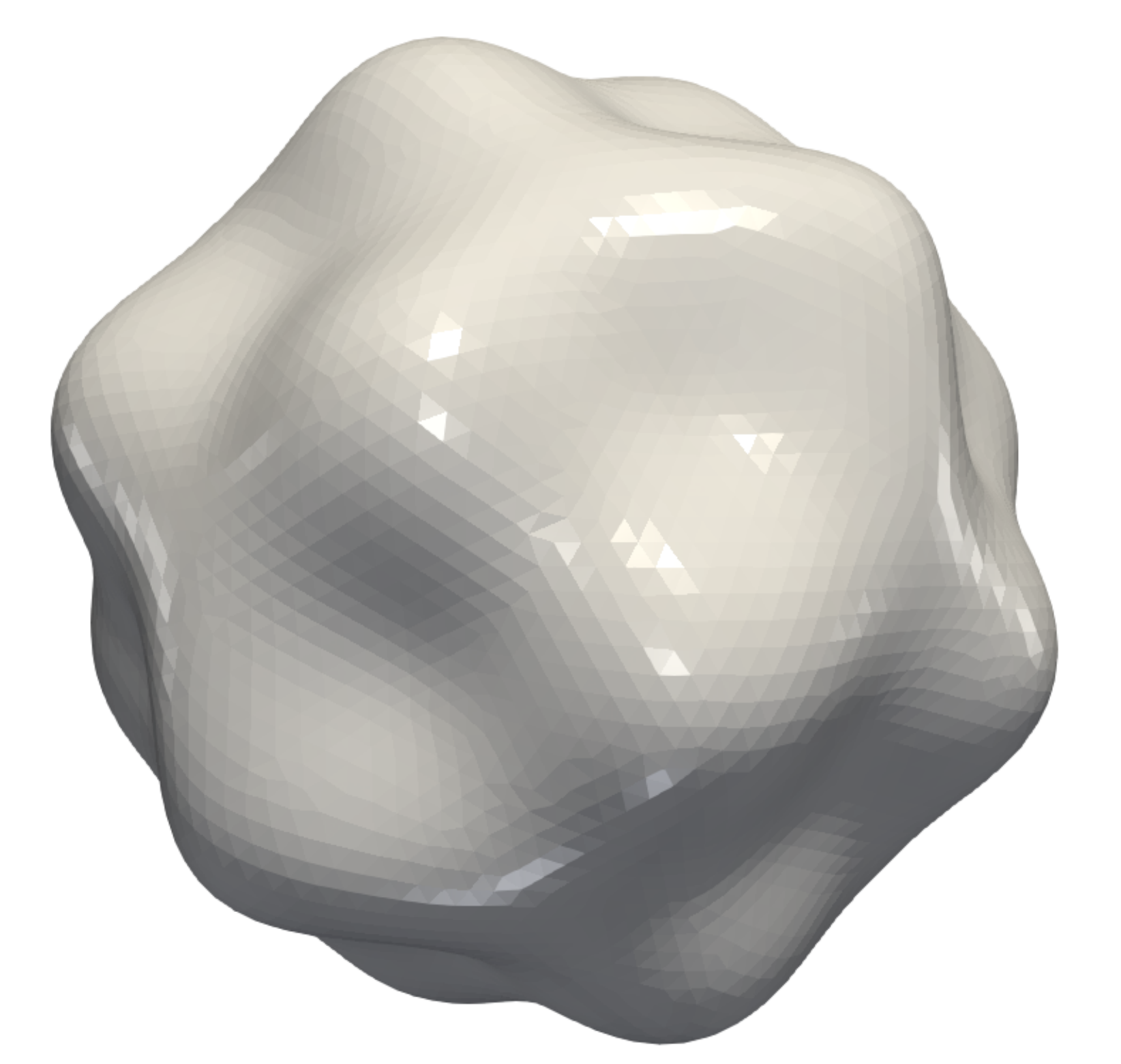}
        \caption{Popcorn surface $\Gamma$.}
        \label{fig:bdry_LM_3d_geometry}
    \end{subfigure}\hfill
    \begin{subfigure}[t]{0.275\linewidth}
        \centering
        \includegraphics[width=\linewidth]{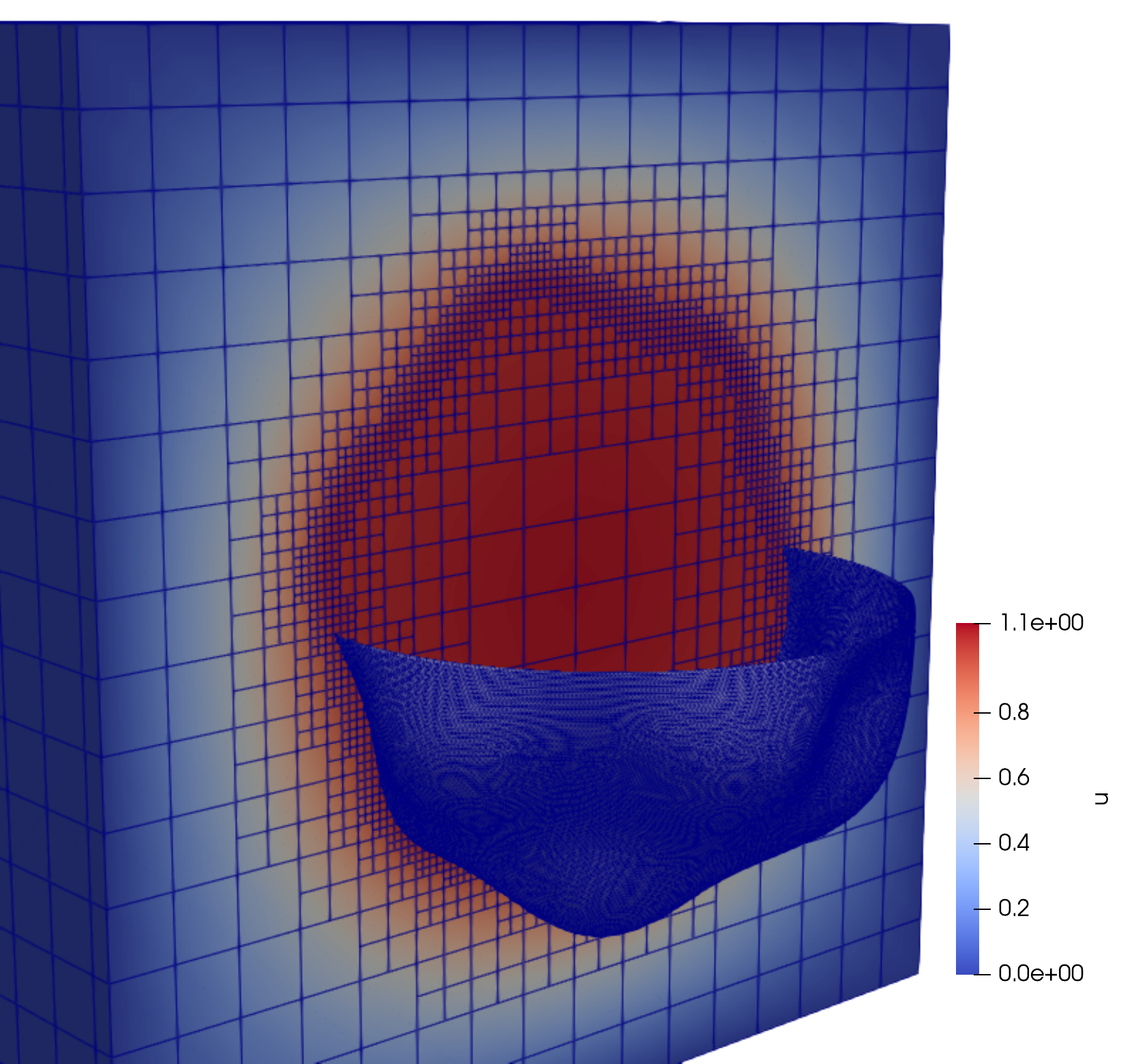}
        \caption{Solution (clipped).}
        \label{fig:bdry_LM_3d_solution}
    \end{subfigure}
    \caption{%
        (\subref{fig:weak_bc_3d_geometry_shell})--(\subref{fig:weak_bc_3d_geometry_sol}): tetrahedralization of the surface mesh (clipped) and computed solution for the weak imposition of boundary conditions with $g=x+y$.
        (\subref{fig:bdry_LM_3d_geometry})--(\subref{fig:bdry_LM_3d_solution}): popcorn surface $\Gamma$ (12,178 triangles) and computed solution on $\Omega=(-1,1)^3$ (clipped in order to show mesh of $\Gamma$) for the Poisson fictitious domain problem with boundary supported Lagrange multipliers, with $f=g=1$, $s=\frac{1}{2}$. Iteration counts are reported in
        Table~\ref{tab:bdry_LM_3d_convergence}.%
    }
    \label{fig:weak_bc_3d_mesh}\label{fig:bdry_LM_3d}
\end{figure}

To stress the robustness of the preconditioner, we consider a non-trivial three-dimensional test case involving the
geometry shown in Figure~\ref{fig:weak_bc_3d_geometry_shell}, keeping $g=x+y$ as boundary datum.  For such a geometry, the generation of a hexahedral mesh is a complex task, and we hence consider tetrahedral elements. Hence, $\mathcal{P}^1$--$\mathcal{P}^1$ is the chosen
finite element pair. The meshing procedure starts from an STL file representing the boundary of the domain\footnote{Taken from \url{https://github.com/alecjacobson/common-3d-test-models}.}, and generates a volumetric tetrahedral mesh through the \textsc{ftetwild} library~\cite{fTetWild}. The resulting mesh is shown in Figure~\ref{fig:weak_bc_3d_geometry_shell}, while the computed
solution is shown in Figure~\ref{fig:weak_bc_3d_geometry_sol}.

The results in terms of iteration counts across mesh refinement are collected
in Table~\ref{tab:weak_bc_3d_convergence}, confirming the robustness of the AL preconditioner also for a challenging geometry. The average number of
inner iterations to approximately invert the augmented block mildly increases with mesh refinement, never exceeding nine PCG iterations.

\begin{figure}[htbp]
    \centering
    \captionsetup[table]{font=small,skip=4pt}
    \begin{minipage}[t]{0.49\linewidth}
        \vspace{0pt}
        \centering
        {\small
            \renewcommand{\arraystretch}{1.05}
            \sisetup{group-separator={,}, group-minimum-digits=4}
            \resizebox{\linewidth}{!}{%
                \begin{tabular}{
                        S[table-format=7]
                        S[table-format=4]
                        S[table-format=2]
                        S[table-format=4]
                    }
                    \toprule
                    \multicolumn{4}{c}{\textbf{Weak bcs.\ through multipliers (2D)}}               \\
                    \midrule
                    {$\mathrm{dof}_u$} & {$\mathrm{dof}_\lambda$} & {FGMRES its.} & {Inner (avg.)} \\
                    \midrule
                    81                 & 16                       & 6             & {2}            \\
                    289                & 32                       & 7             & {2}            \\
                    1089               & 64                       & 6             & {2}            \\
                    4225               & 128                      & 6             & {5}            \\
                    16641              & 256                      & 5             & {5}            \\
                    66049              & 512                      & 5             & {5}            \\
                    263169             & 1024                     & 5             & {6}            \\
                    1050625            & 2048                     & 4             & {6}            \\
                    4198401            & 4096                     & 5             & {7}            \\
                    \bottomrule
                \end{tabular}}}
        \captionof{table}{%
            FGMRES(30) iteration counts for the AL preconditioner $\mathcal{\widetilde{P}}_\gamma^{\text{Babu\v{s}ka}}$ applied to
            weakly imposed boundary conditions~\eqref{eqn:LM_discrete_augmented}
            on $\Omega=(0,1)^2$, with $g=x+y$, $\gamma=10$, $s=\frac{1}{2}$.
        }
        \label{tab:weak_bc_2d_convergence}
    \end{minipage}%
    \hfill
    \begin{minipage}[t]{0.49\linewidth}
        \vspace{0pt}
        \centering
        {\small
            \renewcommand{\arraystretch}{1.05}
            \sisetup{group-separator={,}, group-minimum-digits=4}
            \resizebox{\linewidth}{!}{%
                \begin{tabular}{
                        S[table-format=7]
                        S[table-format=6]
                        S[table-format=2]
                        S[table-format=4]
                    }
                    \toprule
                    \multicolumn{4}{c}{\textbf{Weak bcs.\ through multipliers (3D)}}               \\
                    \midrule
                    {$\mathrm{dof}_u$} & {$\mathrm{dof}_\lambda$} & {FGMRES its.} & {Inner (avg.)} \\
                    \midrule
                    5859               & 3887                     & 15            & {6}            \\
                    37356              & 15542                    & 11            & {7}            \\
                    264283             & 62162                    & 11            & {8}            \\
                    1982981            & 248642                   & 8             & {9}            \\
                    \bottomrule
                \end{tabular}}}
        \captionof{table}{%
            FGMRES(30) iteration counts for the AL preconditioner $\mathcal{\widetilde{P}}_\gamma^{\text{Babu\v{s}ka}}$ applied to
            weakly imposed boundary conditions~\eqref{eqn:LM_discrete_augmented}
            on the three-dimensional geometry in Figure~\ref{fig:weak_bc_3d_mesh}, with
            $g=x+y$, $\gamma=10$, $s=\frac{1}{2}$.
        }
        \label{tab:weak_bc_3d_convergence}
    \end{minipage}

\end{figure}

\subsection{Fictitious domain with boundary supported Lagrange multipliers}\label{sec:num_blm}
We consider the saddle-point system associated with the Poisson fictitious domain problem with \textit{boundary} Lagrange multipliers, which we briefly recall from~\cite{Glowinski}.
Let $\omega$ be a closed and bounded domain of $\mathbb{R}^d$, $d=2,3$, with Lipschitz continuous boundary $\Gamma \coloneqq \partial \omega$, and $\Omega \subset \mathbb{R}^d$ a Lipschitz domain such that $\omega$ is compactly contained in $\Omega$. We consider the associated variational formulation from~\cite{Glowinski}, i.e., find $(u,\lambda) \in V \times \Lambda'$
\begin{equation}\label{eqn:AL_prec_BLM}
    \begin{cases}
        (\nabla u,\nabla v)_\Omega+\langle v,\lambda\rangle_{\Lambda\times\Lambda'}
        =(f,v)_\Omega                                 & \forall v\in V=H^1_0(\Omega),                \\
        \langle u,\mu\rangle_{\Lambda\times\Lambda'}
        =\langle g,\mu\rangle_{\Lambda\times\Lambda'} & \forall \mu\in\Lambda'=H^{-\frac12}(\Gamma),
    \end{cases}
\end{equation}for given data $f \in L^2(\Omega)$ and $g \in H^{\frac{1}{2}}(\Gamma)$. Formally, the above problem and its augmented formulation is very close to
the one in the previous section. The main difference lies in the fact that $\Gamma$ is an immersed interface (a surface in the three-dimensional case) not aligned with the background domain $\Omega$, yielding a loss of regularity for $u$, which cannot be expected to be globally in $H^2(\Omega)$~\cite{Glowinski}. Remark~\ref{rmk:exponent} gives again $\W=h_{\Gamma}M$, therefore the blocks of the preconditioner are defined analogously to those in~\eqref{eqn:AL_prec_babuska}. We refer to the preconditioner for this problem as $\widetilde{\mathcal{P}}_{\gamma}^{\mathrm{BLM}}$.

We test it on a challenging three-dimensional immersed geometry. We set $u|_{\partial\Omega}=0$, $f=g=1$, $\Omega=(-1,1)^3$, and take $\Gamma$ to be a surface mesh representing a popcorn-like geometry. The discretization of $\Gamma$ is done by means of a
surface triangular mesh consisting of 12,178 elements, and it is shown in Figure~\ref{fig:bdry_LM_3d_geometry}.
A localized refinement is performed around the interface, and the computed solution is shown in Figure~\ref{fig:bdry_LM_3d_solution}. Iteration counts in Table~\ref{tab:bdry_LM_3d_convergence} remain robust under mesh refinement.
Low and robust iteration counts (not reported here) are also confirmed for a simpler two-dimensional test with $\Omega=(0,1)^2$ and $\Gamma$ given by the boundary of the ball of radius $0.2$ centered at $(0.4,0.4)$.

    {\small \begin{table}[!htbp]
            \centering
            \sisetup{group-separator={,}, group-minimum-digits=4}
            \begin{tabular}{
                    S[table-format=8]
                    S[table-format=6]
                    S[table-format=2]
                    S[table-format=4]
                }
                \toprule
                \multicolumn{4}{c}{\textbf{FD-BLM (3D)}}                                       \\
                \midrule
                {$\mathrm{dof}_u$} & {$\mathrm{dof}_\lambda$} & {FGMRES its.} & {Inner (avg.)} \\
                \midrule
                141179             & 6091                     & 14            & {11}           \\
                577664             & 24358                    & 15            & {11}           \\
                2442031            & 97426                    & 12            & {12}           \\
                10808121           & 389698                   & 12            & {12}           \\
                \bottomrule
            \end{tabular}
            \caption{%
                FGMRES(30) iteration counts for the AL preconditioner $\mathcal{\widetilde{P}}_\gamma^{\mathrm{BLM}}$ applied to the Poisson
                fictitious domain problem with boundary supported Lagrange multipliers. We set $\Omega=(-1,1)^3$, $\Gamma$ the pop-corn geometry, $f=g=1$,
                $\gamma=10$, $s=\frac{1}{2}$, and we use $\mathcal{Q}^1$--$\mathcal{P}^1$ elements.%
            }
            \label{tab:bdry_LM_3d_convergence}
        \end{table}}

\subsection{Fictitious domain with distributed Lagrange multipliers}\label{sec:num_dlm}
We consider the following \emph{elliptic interface problem} with a jump in the coefficients.

Let $\Omega$ be a domain in $\mathbb{R}^d$, where $d \in \{2,3\}$, with a bounded Lipschitz boundary $\partial \Omega$.
Let $\Omega_1$ and $\Omega_2$ be two subdomains of $\Omega$ such that $\overline{\Omega} = \overline{\Omega}_1 \cup \overline{\Omega}_2$, and let the interface $\Gamma = \overline{\Omega}_1 \cap \overline{\Omega}_2$ be Lipschitz continuous. We assume $\overline{\Gamma} \cap \partial \Omega = \emptyset$.

For $i=1,2$, given forcing terms $f_i \colon \Omega_i \rightarrow \mathbb{R}$, and positive coefficients $\beta_i \in L^{\infty}(\Omega_i)$, find $u_1 \colon \Omega_1 \rightarrow \mathbb{R}$ and $u_2 \colon \Omega_2 \rightarrow \mathbb{R}$ such that:
\begin{equation}
    \begin{cases}\label{eqn:model_problem}
        -\operatorname{div} \bigl( \beta_i \nabla u_i \bigr) = f_i
        \>\> \>\>\>\>\>\> \>\>\>\>\>\> \>\>\>\>\>\> \>\>\>\>\>\> \mathrm{in} \> \Omega_i,                                            \\
        \>\> \>\>\>\>\>\>\>\>\>\>\>\>\>\>\>\>\>\>\>\>u_1 = u_2 \>\> \quad\>\>\>\>\>\>\>\>\>\>\>\>\>\>\>\>\>\> \mathrm{on} \> \Gamma, \\
        \>\>\>\>\>\>\beta_1 \nabla u_1 \cdot \mathbf{n}_1 = -\beta_2 \nabla u_2 \cdot \mathbf{n}_2 \>\>\>\>\mathrm{on}\> \Gamma,     \\
        \qquad \qquad  \> \>\>\>u_1 = 0  \>\>\>\>\>\>\>\>\>\>\>\>\>\>\>\>\>\>\>\>\>\> \quad \mathrm{on} \> \partial \Omega ,
    \end{cases}
\end{equation}
where $\mathbf{n}_i$ ($i=1,2$) denotes the unit vector normal to $\Gamma$, pointing out of $\Omega_i$.

Following the fictitious domain discretization developed in~\cite{AURICCHIO201536,BOFFI20142151}, let $V= H_0^1(\Omega)$, $V_2= H^1(\Omega_2)$, and $\Lambda'= [H^1(\Omega_2)]'$.
Given $f\in L^2(\Omega)$, $f_2\in L^2(\Omega_2)$ and $\beta_2 > \beta >0$, let $\widetilde{V}_h\coloneqq V_h\times V_{2,h}$ denote the
product space, whose elements are $\widetilde{v}_h=(v_h,v_{2,h})$.
The formulation requires finding $(\widetilde{u}_h,\lambda_h) = \big((u_h,u_{2,h}),\lambda_h\big)\in \widetilde{V}_h \times\Lambda_h$ such that
\begin{equation}\label{eqn:FDDLM_discrete_augmented}
    \begin{cases}
        a_e(\widetilde{u}_h,\widetilde{v}_h) + (\lambda_h,\, v_h|_{\Omega_2}-v_{2,h})_{\Omega_2} = (f,v_h)_\Omega + (f_2-f,v_{2,h})_{\Omega_2} & \forall \widetilde{v}_h \in \widetilde{V}_h, \\
        (\mu_h,\, u_h|_{\Omega_2}-u_{2,h})_{\Omega_2} = 0                                                                                      & \forall \mu_h\in\Lambda_h,
    \end{cases}
\end{equation}
where  $a_e \big((u_h,u_{2,h}),(v_h,v_{2,h})\big) := (\beta\nabla u_h,\nabla v_h)_\Omega + \big((\beta_2-\beta)\nabla u_{2,h},\nabla v_{2,h}\big)_{\Omega_2}.$

Hence, for this problem it follows $L(u,u_2)\coloneqq u_{|\Omega_2}-u_2$, which weakly enforces $u$ to match $u_2$ on the immersed domain $\Omega_2$.

We denote by $h_\Omega$ the characteristic mesh size of the background mesh on $\Omega$ and by $h_2$ that of the immersed mesh on $\Omega_2$. We assume that both meshes are shape regular and quasi-uniform and satisfy $h_2 \simeq h_\Omega$. Thus, the analogue of the discrete augmented formulation~\eqref{eqn:problem_general_discrete_augmented} is obtained by adding the term $\gamma h_2^{-\alpha}\big(u_h|_{\Omega_2}-u_{2,h},v_h|_{\Omega_2}-v_{2,h}\big)_{\Omega_2}$
to the left-hand side of the first equation. Moreover, following the recipe outlined in Remark~\ref{rmk:exponent}, we set $\W = h_{2}^{2}M$.
In this case, we consider instead the projected augmentation
$$
    \gamma h_2^{-2}\big(\pi_H(u_h|_{\Omega_2}-u_{2,h}),\pi_H(v_h|_{\Omega_2}-v_{2,h})\big)_{\Omega_2},
$$
as an example of the possible alternative mentioned in Section~\ref{sec:motivation}.
By Lemma~\ref{lemma:algebraic_term}, this term is equivalent to the algebraic augmentation and therefore does not alter the discrete solution. Using finite element spaces for which the multiplier space is discontinuous, i.e.  $\mathcal{Q}^1$--$(\mathcal{Q}^1+\mathfrak{B})$--$\mathcal{Q}^0$, where $\mathfrak{B}$ denotes the bubble enrichment (inf-sup stability for this element has been shown in~\cite{BoffiGastaldi2024}), the projected term can be assembled without forming $C^TM^{-1}C$. Indeed, $\pi_H$ reduces to the local projection onto piecewise constants, which is a local operation consisting of averaging the values of $u_h$ over each cell of the immersed mesh $\Omega_2$.

Introducing basis functions $\{\varphi_i\}_{i=1}^n$, $\{\phi_j\}_{j=1}^{l}$, and $\{\psi_k\}_{k=1}^m$ for the finite element spaces $V_h$, $V_{2,h}$, and $\Lambda_h$, respectively, the resulting
matrices are defined as follows:
\begin{equation}\label{eqn:matrices_elliptic_interface}
    \begin{gathered}
        \begin{aligned}
            \widetilde{A}(\gamma)_{ij}     & = (\beta \nabla\varphi_i , \nabla\varphi_j)_{\Omega} + \gamma h_{2}^{-2} ( \pi_H \varphi_i|_{\Omega_2}, \pi_H \varphi_j|_{\Omega_2})_{\Omega_2}, \\
            \widetilde{A}_{2}(\gamma)_{ij} & = ((\beta_2 -\beta) \nabla\phi_i ,\nabla\phi_j)_{\Omega_2} + \gamma h_{2}^{-2} ( \pi_H \phi_i, \pi_H \phi_j)_{\Omega_2},                         \\
            D(\gamma)_{ij}                 & = \gamma h_2^{-2}(\pi_H \varphi_i|_{\Omega_2}, \pi_H \phi_j)_{\Omega_2},                                                                         \\
            C_{ki}                         & = (\varphi_i|_{\Omega_2},\psi_k)_{\Omega_2},                                                                                                     \\
            C_{2_{kj}}                     & = (\phi_j,\psi_k)_{\Omega_2},                                                                                                                    \\
            M_{ij}                         & = (\psi_i,\psi_j)_{\Omega_2}.
        \end{aligned}
    \end{gathered}
\end{equation}

Our framework gives the following candidate (ideal) preconditioner, where no triple matrix product needs to be formed:
\begin{equation}\label{eqn:ideal_AL_elliptic_interface_matrix_form}
    \mathcal{\widetilde{P}}_{\gamma}^{\mathrm{DLM}}\coloneqq
    \begin{pmatrix}
        \widetilde{A}(\gamma) & - D(\gamma)             & C^T                  \\
        - D(\gamma)^T         & \widetilde{A}_2(\gamma) & -C_2^T               \\
        0                     & 0                       & -\frac{1}{\gamma} \W
    \end{pmatrix}.
\end{equation}

Nevertheless, approximately inverting the augmented $2\times2$ block remains non-trivial. To address this issue, the authors introduced in~\cite{benzi2026augmentedlagrangianpreconditionersfictitious} a \emph{modified} version with two distinct augmentation parameters, $\gamma_1$ and $\gamma_2$. In our framework, this preconditioner takes the form
\begin{equation}\label{eqn:modified_AL_elliptic_interface_matrix_form}
    \mathcal{\widetilde P}_{\gamma_1,\gamma_2}^{\mathrm{DLM}}\coloneqq
    \begin{pmatrix}
        \widetilde{A}(\gamma_1) & - D(\gamma_1)             & C^T                    \\
        0                       & \widetilde{A}_2(\gamma_2) & -C_2^T                 \\
        0                       & 0                         & -\frac{1}{\gamma_1} \W
    \end{pmatrix},
\end{equation}and only requires the inversion of the SPD diagonal blocks. As shown in~\cite{benzi2026augmentedlagrangianpreconditionersfictitious}, the preconditioner
is fairly robust with respect to the jump $\beta_2 - \beta$, with an increasingly tighter spectrum when $\gamma_1 \rightarrow \infty$ and $\gamma_2 \rightarrow 0$. For this reason, following~\cite{benzi2026augmentedlagrangianpreconditionersfictitious}, we
set $\gamma_1=10$ and $\gamma_2= 10^{-2}$.

\begin{figure}[htbp]
    \centering
    \begin{minipage}[c]{0.55\linewidth}
        \centering
        \sisetup{group-separator={,}, group-minimum-digits=4}
        \resizebox{\linewidth}{!}{%
            \begin{tabular}{
                    S[table-format=7]
                    S[table-format=7]
                    S[table-format=6]
                    S[table-format=2]
                    S[table-format=2]
                    S[table-format=2]
                    S[table-format=4]
                }
                \toprule
                \multicolumn{7}{c}{\textbf{FD-DLM (2D)}}                                                                                                                   \\
                \midrule
                {$\mathrm{dof}_u$} & {$\mathrm{dof}_{u_2}$} & {$\mathrm{dof}_\lambda$} & \multicolumn{3}{c}{FGMRES iterations} & {Inner (avg.)}                            \\
                \cmidrule(lr){4-6}
                                   &                        &                          & {$\beta_2=10$}                        & {$\beta_2=10^3$} & {$\beta_2=10^7$} &     \\
                \midrule
                4225               & 169                    & 80                       & 18                                    & 18               & 18               & {6} \\
                16641              & 657                    & 320                      & 23                                    & 23               & 23               & {6} \\
                66049              & 2593                   & 1280                     & 22                                    & 22               & 22               & {7} \\
                263169             & 10305                  & 5120                     & 22                                    & 22               & 22               & {7} \\
                1050625            & 41089                  & 20480                    & 23                                    & 23               & 23               & {7} \\
                4198401            & 164097                 & 81920                    & 25                                    & 24               & 24               & {8} \\
                \bottomrule
            \end{tabular}}%
        \captionof{table}{%
            FGMRES(30) iteration counts for the \emph{modified} AL preconditioner $\mathcal{\widetilde{P}}_{\gamma_1,\gamma_2}^{\mathrm{DLM}}$ applied to the elliptic interface
            fictitious domain problem with distributed Lagrange multipliers. We set
            $\Omega=(-1,1)^2$, $\Omega_2=\mathcal{B}_{0.3}(0,0)$,
            $\mathcal{Q}^1$--$(\mathcal{Q}^1 + \mathfrak{B})$--$\mathcal{Q}^0$ elements, for increasing values of
            $\beta_2$. The augmentation parameters are chosen as $\gamma_1 = 10$, $\gamma_2 = 10^{-2}$, $s=1$.
        }
        \label{tab:dlm_2d_convergence_q0}
    \end{minipage}%
    \hfill
    \begin{minipage}[c]{0.40\linewidth}
        \centering
        \includegraphics[width=\linewidth]{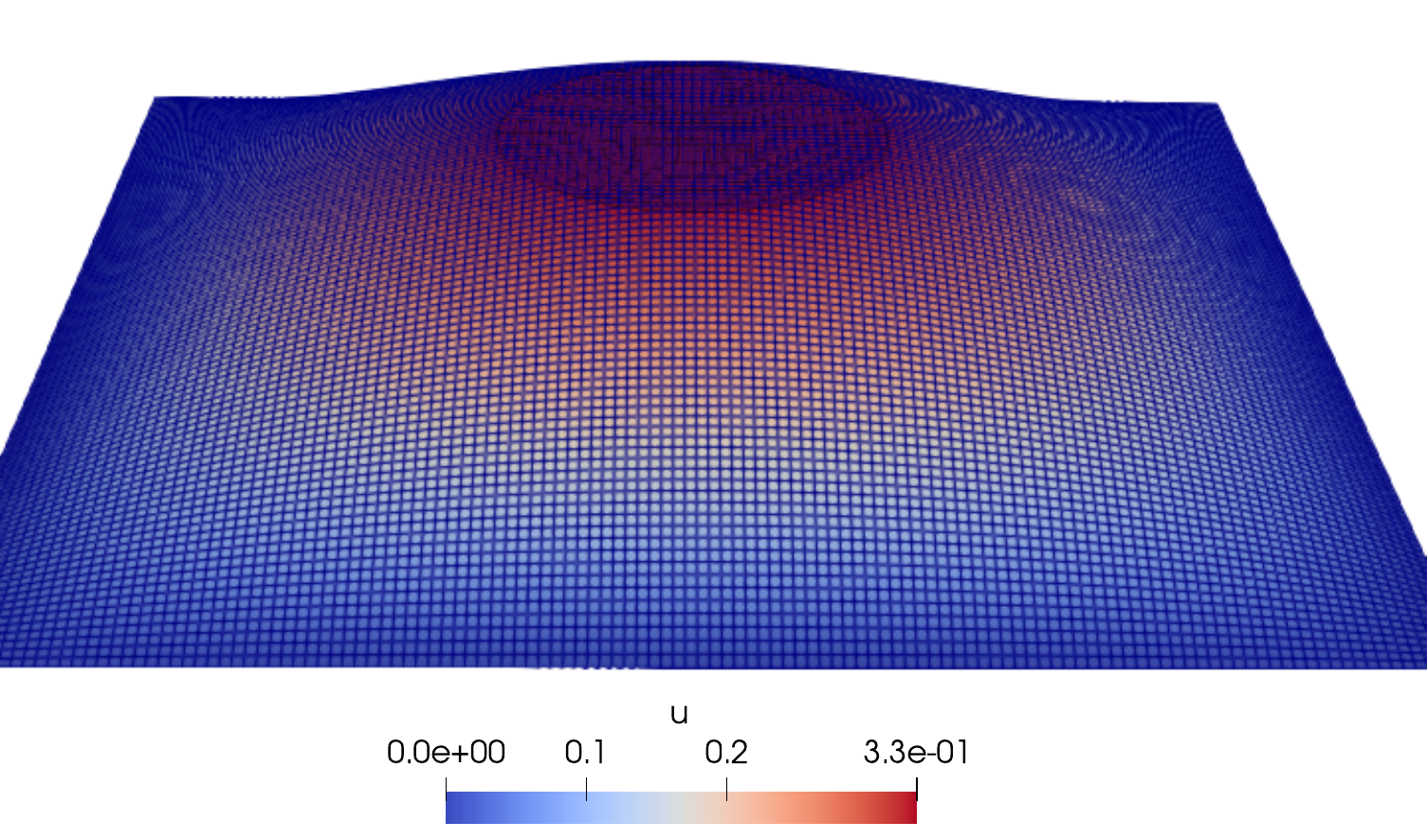}
        \caption{%
            Computed solution $u_h$ for the elliptic interface fictitious domain problem with distributed
            Lagrange multipliers ($\beta_2\!=\!10^3$) on $\Omega=(-1,1)^2$, $\Omega_2=\mathcal{B}_{0.3}(0,0)$.
        }
        \label{fig:dlm_2d}
    \end{minipage}
\end{figure}

To validate this preconditioner, we start with a classic two-dimensional experiment where $\Omega=(-1,1)^2$ and $\Omega_2 = \mathcal{B}_{0.3}(0,0)$, the ball of radius $0.3$ centered at the origin. The number of FGMRES iterations
for different values of $\beta_2$ are collected in Table~\ref{tab:dlm_2d_convergence_q0}, showing high robustness of the preconditioner also with respect to the jump in the coefficients induced by parameter $\beta_2$.
Inner iteration counts for the approximate inversion of the augmented (1,1)-block when $\beta_2=10^3$ are reported in the last column and low values are observed across mesh refinement.

A visualization of the solution for $\beta_2=10^3$ is shown in Figure~\ref{fig:dlm_2d}, where the contrast among subdomains results in a clear and visible separation of the two regions.

For the three-dimensional case, we consider a more complex setting by extending the scalar elliptic interface problem to an immersed linear elasticity problem. We use the same finite element spaces as in the previous scalar case.
Let $\mu,\mu_2$ be the shear moduli in $\Omega$ and $\Omega_2$, respectively, and $\lambda, \lambda_2$ the first Lamé parameters. Consistently with the assumptions for the scalar problem~\eqref{eqn:model_problem}, we assume $\mu_2 > \mu>0$ and $\lambda_2 > \lambda>0$. We denote with $\delta_\mu = \mu_2 - \mu$ and $\delta_\lambda = \lambda_2 - \lambda$ the corresponding jumps. Vector-valued functions are indicated in boldface; thus, no ambiguity arises when $\lambda$ denotes either a multiplier or a Lamé parameter. Adapting the arguments from the previous sections
to the vector-valued case, the fictitious domain formulation reads as follows: find $(\boldsymbol{u},\boldsymbol{u}_2,\boldsymbol{\lambda})\in V\times V_2\times\Lambda'$ such that
\begin{align*}
    2 \mu(\varepsilon(\boldsymbol{u}),\varepsilon(\boldsymbol{v}))_\Omega + \lambda(\operatorname{div} \boldsymbol{u}, \operatorname{div} \boldsymbol{v})_\Omega + \langle \boldsymbol{v}|_{\Omega_2}, \boldsymbol{\lambda} \rangle_{\Lambda\times\Lambda'}                        & =(\boldsymbol{f},\boldsymbol{v})_\Omega &  & \forall \boldsymbol{v} \in V,        \\
    2 \delta_\mu (\varepsilon(\boldsymbol{u}_2),\varepsilon(\boldsymbol{v}_2))_{\Omega_2} + \delta_\lambda(\operatorname{div} \boldsymbol{u}_2, \operatorname{div} \boldsymbol{v}_2)_{\Omega_2} - \langle \boldsymbol{v}_{2}, \boldsymbol{\lambda} \rangle_{\Lambda\times\Lambda'} & = 0                                     &  & \forall \boldsymbol{v}_{2}\in V_{2}, \\
    \langle \boldsymbol{u}|_{\Omega_2} - \boldsymbol{u}_{2}, \boldsymbol{\mu} \rangle_{\Lambda\times\Lambda'}                                                                                                                                                                      & = 0                                     &  & \forall \boldsymbol{\mu}\in\Lambda',
\end{align*}

where $\varepsilon(\boldsymbol{u}) = \frac{1}{2}(\nabla \boldsymbol{u} + \nabla \boldsymbol{u}^{\mathsf{T}})$ is the linearized strain tensor, $\boldsymbol{f} \in \bigl[L^2(\Omega)\bigr]^3$, and $V = [H_{\Gamma_\mathrm{D}}^1(\Omega)]^3, V_2=[H^1(\Omega_2)]^3, \Lambda'=\bigl[[H^{1}(\Omega_2)]^{3}\bigr]'$ are the vector-valued analogues of the spaces defined for the scalar problem.
Here, $\Gamma_\mathrm{D}$ denotes the portion of the boundary on which Dirichlet conditions are imposed.

\begin{figure}[tbp]
    \centering
    \begin{subfigure}[b]{0.42\linewidth}
        \centering
        \includegraphics[height=4cm,keepaspectratio]{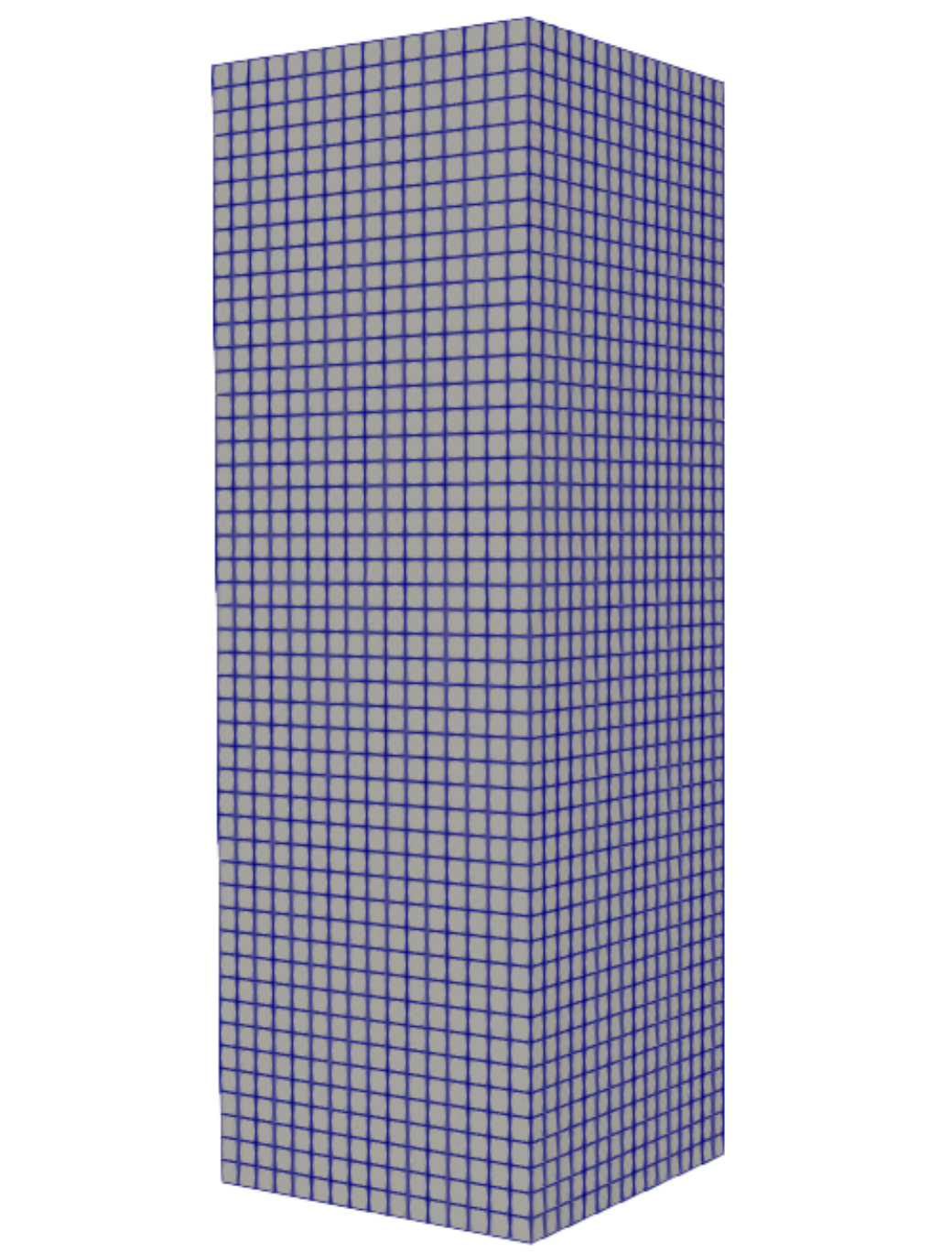}
        \caption{Rectangular inclusion for $\Omega_2$.}
        \label{fig:dlm_3d_geometry}
    \end{subfigure}
    \hspace{0.04\textwidth}
    \begin{subfigure}[b]{0.42\linewidth}
        \centering
        \includegraphics[height=4cm,keepaspectratio]{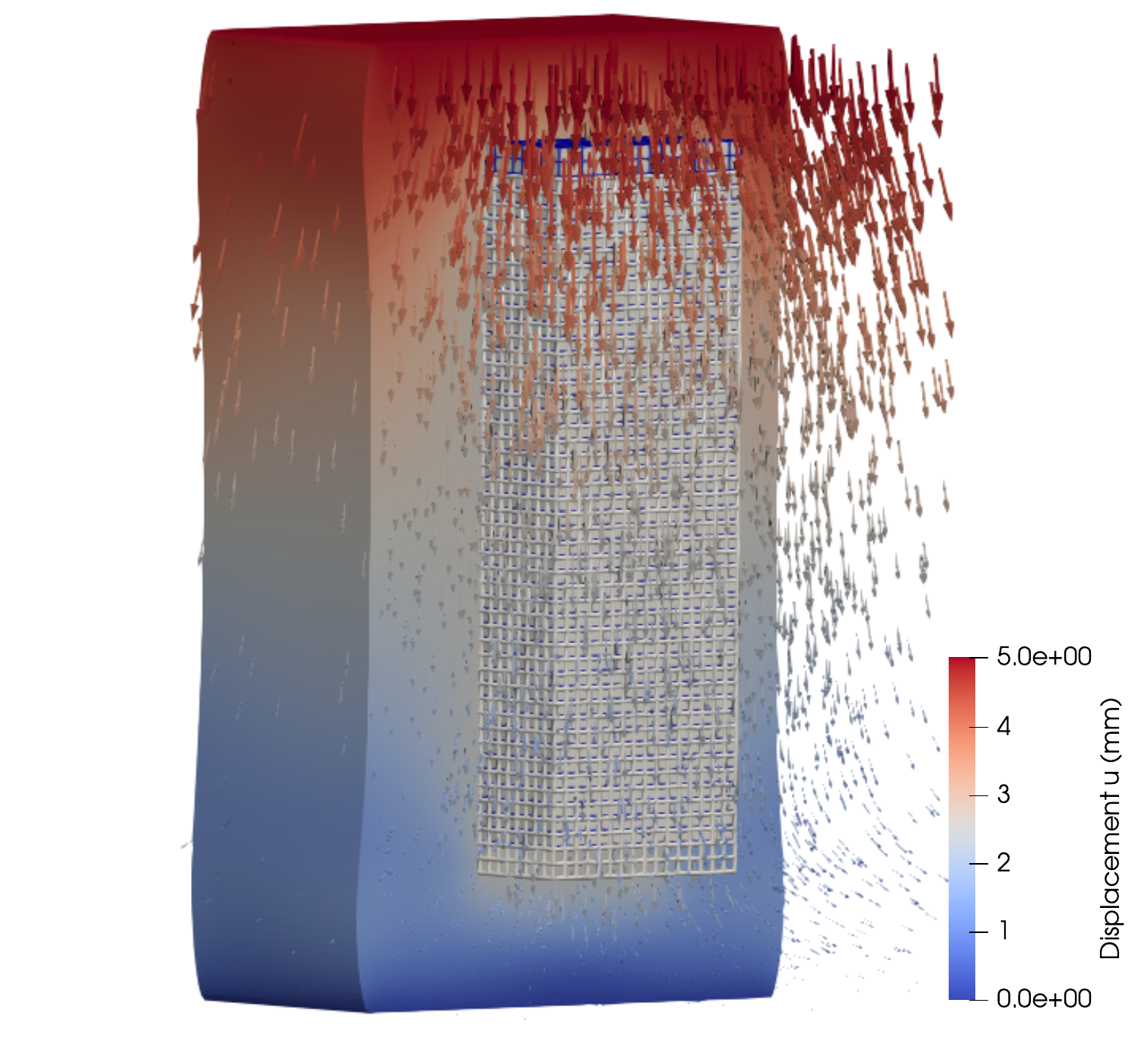}
        \caption{Computed solution (clipped).}
        \label{fig:dlm_3d_solution}
    \end{subfigure}
    \caption{%
        Linear elasticity problem: bar inclusion mesh immersed in a rubber-like elastomeric
        mount. In gray, the displaced configuration of the inclusion, while in blue the reference configuration. The background domain
        is deformed and clipped in order to show the immersed geometry.
    }
    \label{fig:dlm_3d}
\end{figure}

We consider an elongated steel bar inside a rubber-like elastomeric mount.
The background domain is the $\SI{40}{\milli\metre} \times
    \SI{40}{\milli\metre} \times \SI{80}{\milli\metre}$ rectangular box
$\Omega = (-20,20) \times (-20,20) \times (-40,40)$ (in \si{\milli\metre}), while
the immersed domain is the $\SI{18}{\milli\metre} \times \SI{18}{\milli\metre}
    \times \SI{54}{\milli\metre}$ box $\Omega_2 = (-9,9) \times (-9,9) \times (-27,27)$. Both domains are
discretized with structured hexahedral meshes, with the bar mesh shown in Figure~\ref{fig:dlm_3d_geometry}.

The material properties are chosen to be representative of the physical materials: the bar has Lamé constants $\mu_2=\SI{76923}{\mega\pascal}$ and $\lambda_2= \SI{115385}{\mega\pascal}$, corresponding to a Young's modulus $E_2 = \SI{200}{\giga\pascal}$, while the background domain has $\mu = \SI{17.24}{\mega\pascal}$ and $\lambda = \SI{155.17}{\mega\pascal}$, corresponding to a Young's modulus $E = \SI{50}{\mega\pascal}$. This implies
a contrast of nearly four orders of magnitude between the immersed inclusion and the background domain, and represents a good
test case to evaluate the AL preconditioner on a physically relevant scenario. On the background domain $\Omega$, we prescribe a fixed displacement field $\boldsymbol{u} = (0,0,-\delta)$, with $\delta = \SI{5}{\milli\metre}$, on the top face and zero displacement on the bottom face. Traction-free boundary conditions ($\boldsymbol{\sigma} \cdot \boldsymbol{n} = 0$)
are applied to the four remaining lateral faces, and body forces are neglected. The displacement magnitude (in \si{\milli\metre}) with the associated streamlines is shown
in Figure~\ref{fig:dlm_3d_solution}. The softer background domain shows a smooth gradient for the magnitude of the displacement from the clamped bottom to the fixed top face, typical of a uniaxial compression state. Due to the high
contrast in material properties, the bar undergoes a translation of $\SI{2.5}{\milli\metre}$ in the negative $z$-direction with no visible distortion, behaving as a rigid body.

We use the same inner PCG solvers as in the previous experiments, with an outer relative residual tolerance of $10^{-6}$. We then study the behavior of $\mathcal{\widetilde P}_{\gamma_1,\gamma_2}^{\mathrm{DLM}}$ with respect to simultaneous mesh refinement up to $13 \times 10^6$
degrees of freedom. The FGMRES iterations for different cycles
are reported in Table~\ref{tab:dlm_3d_convergence}: also in this complex scenario, the iteration counts remain low and robust under mesh refinement. The low number of iterations, combined with a realistic setup featuring high contrast in the material coefficients, motivates the
usage of the proposed preconditioner as a viable option for the solution of large-scale three-dimensional problems.

\begin{table}[!ht]
    \centering
    \sisetup{group-separator={,}, group-minimum-digits=4}
    \begingroup
    \small
    \setlength{\tabcolsep}{4pt}
    \begin{tabular}{
            S[table-format=8]
            S[table-format=7]
            S[table-format=8]
            S[table-format=2]
            S[table-format=4]
        }
        \toprule
        \multicolumn{5}{c}{\textbf{FD-DLM (3D)}}                                                                \\
        \midrule
        {$\mathrm{dof}_u$} & {$\mathrm{dof}_{u_2}$} & {$\mathrm{dof}_\lambda$} & {FGMRES its.} & {Inner (avg.)} \\
        \midrule
        28611              & 1551                   & 576                      & 12            & {7}            \\
        212355             & 10683                  & 4608                     & 10            & {7}            \\
        1635075            & 79347                  & 36864                    & 12            & {9}            \\
        12830211           & 611811                 & 294912                   & 9             & {9}            \\
        \bottomrule
    \end{tabular}
    \endgroup
    \caption{%
        FGMRES(30) iteration counts for the modified AL preconditioner $\mathcal{\widetilde{P}}_{\gamma_1,\gamma_2}^{\mathrm{DLM}}$ applied to the linear elasticity
        fictitious domain problem with distributed Lagrange multipliers on the rectangular
        inclusion immersed in a rubber-like elastomeric mount.
    }
    \label{tab:dlm_3d_convergence}
\end{table}

\subsection{Strong scaling study}
We conclude with a strong scaling experiment, using the same inner and outer solvers and stopping criteria as in the previous experiment. The simulations have been run on the MareNostrum 5 cluster at Barcelona Supercomputing Center, where each compute node features 2x Intel Sapphire Rapids 8460Y+ processors at 2GHz (56 cores each, 112 cores per node) and 256GB of main memory\footnote{Technical characteristics are available at \url{https://www.bsc.es/marenostrum/marenostrum-5}.}.

We consider the same material parameters and spaces as in the previous example, with $\Omega=(-1,1)^3$ and $\Omega_2=(-0.8,0.8)\times(-0.3,0.3)\times(-0.4,0.5)$. The discretization has $50,923,779$ dof for $V_h$, $206,115$ for $V_{2,h}$, and $98,304$ for $\Lambda_h$. We vary the number of MPI processes from 336 to 5376 (i.e., from 3 to 48 compute nodes) and report the wall-clock time as a function of the number of MPI processes for each phase: linear solver, assembly of background and immersed matrices, and coupling terms.

\begin{figure}[!htbp]
    \centering
    \includegraphics{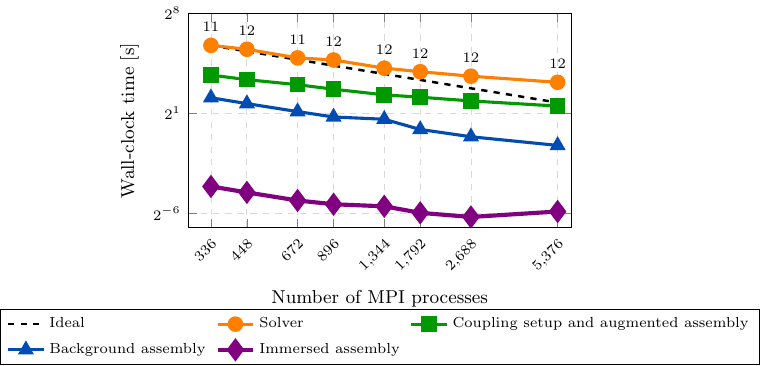}
    \caption{%
        Strong scaling analysis for the fictitious domain method with distributed Lagrange
        multipliers on MareNostrum~5. The dashed line represents ideal scaling. Black numbers above
        the solver curve indicate FGMRES iteration counts.
    }
    \label{fig:strong_scaling}
\end{figure}

Results are shown in Figure~\ref{fig:strong_scaling}. The annotations (in black) indicate the number of FGMRES iterations, which is independent of the number of MPI processes. The linear solver and the classical assembly
phases exhibit good scaling, with a degradation when communication overhead becomes significant due to the low number of dof per MPI process. Particle detection and augmented term assembly dominate the assembly cost. The detection phase locates immersed-mesh quadrature points, treated as particles, in the background mesh and assigns them to subdomains, requiring substantial global communication between the independently partitioned meshes.

\section{Conclusions, practical guidelines and outlook}\label{sec:conc}
We have introduced a unified framework for the design and analysis of augmented Lagrangian preconditioners for symmetric saddle point systems arising from mixed finite element discretizations. The proposed approach is based on a continuous augmentation, which combines the spectral benefits of algebraic augmentation with the computational simplicity of a standard finite element formulation. The augmented $(1,1)$-block can be explicitly assembled, preserving the sparsity pattern of the original block and avoiding sparse matrix-matrix products. Moreover, if we can prove that the stabilizing term vanishes under mesh refinement, then the solution of the discrete augmented formulation converges to that of the original problem. In addition, our analysis shows that the matrix preconditioned with the ideal preconditioner has a real, positive spectrum that is bounded away from zero and clusters around 1 as the AL parameter increases.

Numerical experiments with mixed-dimensional Lagrange multiplier formulations, different discretizations, and complex geometries validate the theoretical findings, while the abstract analysis applies to a broader class of constrained formulations.
The applicability of the framework to a new problem can be assessed through four main checks: coercivity~\eqref{eqn:coercivity_augmented_eq}, the inverse inequality~\eqref{eqn:inequality_multiplier}, the vanishing of the augmentation residual in~\eqref{eqn:decay_stabilizing_term} under mesh refinement, and compatibility of the mesh-dependent scaling with the convergence rate (Section~\ref{subsubsec:weak_bc_alpha_power}).
These conditions provide a practical guideline for assessing whether the framework provides an efficient alternative to algebraic augmentation while retaining its spectral advantages.
The procedure is illustrated in Section~\ref{subsec:weak_bcs} for the weak imposition of boundary conditions considered.

If one of the previous checks does not hold, a possible alternative is to employ the projected counterpart of the augmentation term, as illustrated in Section~\ref{sec:num_dlm}. This approach still avoids the explicit computation of the triple matrix product $C^T W^{-1}C$.

Future work will extend this approach to fluid-structure interaction problems and develop tailored multigrid solvers for the approximate inversion of the augmented blocks.

\section*{Acknowledgments}
The authors are grateful to Patrick Farrell (University of Oxford) for helpful discussions.

\section*{Declarations}
The implementation of this methodology is available in the authors' maintained GitHub repository: \url{https://github.com/fdrmrc/fictitious_domain_AL_preconditioners}. The authors used
ChatGPT v6.1 to improve readability and language. They reviewed and edited the material. The authors take full responsibility for all content.

\bibliographystyle{siamplain}
\bibliography{refs}
\end{document}